\documentclass[11pt]{amsart}
\usepackage[top=3.5cm, bottom=3.5cm, left=2.5cm, right=2.5cm]{geometry}

\usepackage[utf8]{inputenc}
\usepackage{amsmath}
\usepackage{amsfonts}
\usepackage{amssymb}
\usepackage{amsthm}
\usepackage{dsfont}
\usepackage{hyperref}
\usepackage[capitalise,compress]{cleveref}
\usepackage{tikz}
\usepackage{enumitem}

\numberwithin{equation}{section}

\theoremstyle{plain}
\newtheorem{theorem}{Theorem}[section]
\newtheorem{lemma}[theorem]{Lemma}
\newtheorem{proposition}[theorem]{Proposition}
\newtheorem{corollary}[theorem]{Corollary}

\theoremstyle{definition}
\newtheorem{definition}{Definition}[section]

\newtheorem*{restate}{Theorem 2.3}

\theoremstyle{remark}
\newtheorem*{remark}{Remark}
\newtheorem*{note}{Remark}
\newtheorem{Remark}[theorem]{Remark}

\newcommand{\N}{\mathbb{N}}
\newcommand{\Z}{\mathbb{Z}}
\newcommand{\R}{\mathbb{R}}

\newcommand{\Prob}{\mathbb{P}}
\newcommand{\id}{\mathrm{id}}
\newcommand{\eps}{\varepsilon}

\newcommand{\Aut}{\mathrm{Aut}}

\renewcommand{\ge}{\geqslant}
\renewcommand{\le}{\leqslant}
\renewcommand{\geq}{\geqslant}
\renewcommand{\leq}{\leqslant}

\title{Spread-out percolation on transitive graphs of polynomial growth}
\author{Panagiotis Spanos}
\address{Faculty of Mathematics, Ruhr University Bochum, Germany}
\email{panagiotis.spanos@ruhr-uni-bochum.de}

\author{Matthew Tointon}
\address{School of Mathematics, University of Bristol, United Kingdom}
\email{m.tointon@bristol.ac.uk}

\date{}
\renewcommand\footnotemark{}

\begin{document}
\thanks{The first author acknowledges the support of the Austrian Science Fund (FWF): W1230. This work was also partially supported by a Small Grant from the Heilbronn Institute for Mathematical Research and by the German Research Foundation (DFG): SPP2265.}
\maketitle

\begin{abstract}
Let $G$ be a vertex-transitive graph of superlinear polynomial growth. Given $r>0$, let $G_r$ be the graph on the same vertex set as $G$, with two vertices joined by an edge if and only if they are at graph distance at most $r$ apart in $G$. We show that the critical probability $p_c(G_r)$ for Bernoulli bond percolation on $G_r$ satisfies $p_c(G_r)\sim1/\deg(G_r)$ as $r\to\infty$. This extends work of Penrose and Bollobás--Janson--Riordan, who considered the case $G=\Z^d$.

Our result provides an important ingredient in parallel work of Georgakopoulos in which he introduces a new notion of dimension in groups. It also verifies a special case of a conjecture of Easo and Hutchcroft.
\end{abstract}

\setcounter{tocdepth}{1}
\tableofcontents

\section{Introduction}
It is well known that the \emph{critical percolation probability} $p_c(G)$ of a graph $G$ satisfies
\begin{equation}\label{eq:pc>}
p_c(G)\ge\frac1{\deg(G)-1},
\end{equation}
where $\deg(G)$ denotes the maximum degree of a vertex of $G$. Among \emph{transitive} graphs, exact equality is achieved if and only if $G$ is a tree (see \cref{prop:tree}). However, there are various natural families of transitive graphs with degrees tending to infinity that are very far from being trees and nonetheless achieve equality \emph{asymptotically} in the sense that
\begin{equation}\label{eq:pc.deg.asymp}
p_c(G)=\frac{1+o(1)}{\deg(G)},
\end{equation}
where the error term $o(1)$ tends to zero as the degree tends to infinity. For example, Kesten \cite{kesten90} showed that \eqref{eq:pc.deg.asymp} holds for the hypercubic lattice $\Z^d$ as $d\to\infty$ (see also \cite[\S4]{abs}). Penrose \cite[Theorem 1]{Pen} showed that it holds for the `spread-out' lattice, first considered in this context by Hara and Slade \cite{HaraSlade}, in which $x,y\in\Z^d$ are joined by an edge if and only if $\|x-y\|\le r$, with $d>1$ fixed and $\|\,\cdot\,\|$ some norm on $\R^d$, as $r\to\infty$.

In this paper we extend Penrose's result in a natural way to all transitive graphs of polynomial growth. The class of transitive graphs of polynomial growth includes all Cayley graphs on $\Z^d$, and more generally all Cayley graphs on nilpotent groups, but also includes graphs that are not realizable as Cayley graphs. This result lies within the framework of Benjamini and Schramm's influential paper `Percolation beyond $\mathbb{Z}^d$' \cite{bperc96}. The study of probability in general, and percolation in particular, specifically on transitive graphs of polynomial growth is also well established and has seen significant activity in recent years (see e.g. \cite{GeA2,ItR,Em1,Em2,CMT.supercritical,CMT.locality,Dia,gorski,Sat,Sat2,Sat3}).

We define spread-out percolation on a transitive graph $G=(V,E)$ in a similar way to Penrose, but with the graph distance $d_G$ in place of the norm $\|\,\cdot\,\|$. Precisely, for each $r>0$ we consider the graph $G_r$ with vertex set $V$, and with $u,v\in V$ joined by an edge if and only if $d_G(u,v)\le r$. Given $r\ge0$ and $v\in V$, write $\beta_G(r)=\#\{u\in V:d_G(u,v)\le r\}$, the number of vertices in a ball of radius $r$ in $G$. We say $G$ has \emph{polynomial growth} if there exist $C,d>0$ such that $\beta_G(r)\le Cr^d$ for all $r\in\N$, and \emph{superlinear growth} if $\beta_G(r)/r\to\infty$ as $r\to\infty$. 
\begin{theorem}\label{thm:pc}
Let $G$ be a transitive graph with superlinear polynomial growth. Then
\begin{equation}\label{eq:pc}
p_c(G_r)=\frac{1+o(1)}{\beta_G(r)}
\end{equation}
as $r\to\infty$.
\end{theorem}
The assumption of superlinear growth here is necessary, since $p_c(G_r)=1$ for all $r$ if $G$ does not have superlinear growth. Note that $\deg(G_r)=\beta_G(r)-1$, so that \eqref{eq:pc} is equivalent to \eqref{eq:pc.deg.asymp} for the family $(G_r)_{r\in\N}$.

\medskip

Further motivation for this work arises in connection with a new notion of dimension in groups (and more general metric spaces) called the \emph{percolation dimension}, introduced in parallel work by Georgakopoulos \cite{pdim}. This notion is designed to behave similarly to the growth rate, but to assign a more `realistic' value of dimension in the case of graphs with dead ends, which can allow graphs that are `low dimensional' in one sense to nonetheless exhibit very fast growth. The percolation dimension has a number of desirable properties, such as being monotone decreasing with respect to subgroups and quotients. In groups of polynomial growth, it also agrees with the degree of growth, and it is this last property where our work comes in: the fact that the percolation dimension is bounded below by the degree of polynomial growth for such groups is an immediate corollary of \cref{thm:pc}.

\cref{thm:pc} also provides evidence in favour of a recent conjecture of Easo and Hutchcroft. They show that transitive graphs do not satisfy \eqref{eq:pc.deg.asymp} in general by constructing a family $(G_n)$ of transitive graphs with $\deg(G_n)\to\infty$ such that $p_c(G_n)>3/\deg(G_n)$ for large enough $n$ \cite[Figure~5]{locality}, but conjecture that there should nonetheless exist a universal constant $C$ such that if $G$ is an infinite, connected, transitive, simple graph with superlinear growth then $p_c(G)\le C/\deg(G)$ \cite[Conjecture 7.3]{locality}. \cref{thm:pc} implies that for any given transitive graph $G$ of polynomial growth the familiy $(G_r)_{r\in\N}$ satisfies this conjecture.

Lastly, we wish to point out that the proof offers a systematic approach to studying random processes on graphs with polynomial growth, along with a new result that could potentially be applied in other contexts.
Let us elaborate further. 
Our approach to proving \cref{transitive} involves adapting an alternative proof of Penrose's result given by Bollobás, Janson and Riordan \cite{BJR}. 
In their work, they treat the graph locally using their earlier work on inhomogeneous random graphs \cite{BJR2}.
The percolation configuration is determined by a kernel (defined in Section \ref{4}), from which they derive a local behaviour for the percolated graph, which is the key component of their proof.
This approach relies on a trivial property of the Lebesgue measure: 
that the sequence of discrete measures assigning a mass $1/r^d$ to the elements of $(1/r)\cdot \mathbb{Z}^d$ converges weakly to a scaled Lebesgue measure as $r$ tends to infinity.
We generalize this property in the framework of nilpotent groups.
For a torsion-free nilpotent group $\Gamma$, there exists a growth degree $d_\Gamma$, a Lie group $G_\infty$ in which $\Gamma$ can be embedded as a discrete subset, and dilation functions $\delta_\rho$ for $\rho>0$.
We establish that the sequence of discrete measures assigning mass $1/r^{d_\Gamma}$ to the elements of $\delta_{(1/r)}(\Gamma)$ converges weakly to the Haar measure with a scaling constant that we determine precisely (see \cref{Haartheo}. We believe this result has scope to be used in other settings.

\subsection{Organization of the paper} In Section \ref{2-1} we provide a sketch of the proof and precise statements of the results. In Section \ref{2}, we provide a description of the constructions on nilpotent groups, and known results about the metrics and the asymptotics of those groups.  In Section \ref{3} we prove asymptotic properties of a torsion-free nilpotent group, including the proof of Theorem \ref{Haartheo}. In Section \ref{4} we provide a brief summary of the results we need on percolation and from Bollobás, Janson and Riordan's work on inhomogeneous random graphs. In Section \ref{5} we prove \cref{4.1}, and in Section \ref{6} we deduce our main result, \cref{transitive}.

\subsection*{Acknowledgements}
We thank Agelos Georgakopoulos for introducing us to the works of Penrose and Bollobás--Janson--Riordan, and for pointing out their potential connection to his work on percolation dimension. We also thank Tom Hutchcroft and Wolfgang Woess for useful discussions, Anders Karlsson for helpful suggestions, and Yago Moreno Alonso for corrections to a draft. PS thanks Emmanuel Breuillard for hosting him at the University of Cambridge, and the School of Mathematics at the University of Bristol for their hospitality during two visits.

\section{Main result and setting}\label{2-1}
Our proof of \cref{thm:pc} uses a rather striking array of deep machinery from different areas of mathematics. In this section we give a detailed overview, including references.

Let us first restate \cref{thm:pc} in a slightly different form, more along the lines of Penrose's work and the work of Bollobás, Janson and Riordan we mentioned briefly above. Given a transitive graph $G=(V,E)$ and a parameter $\lambda\ge0$, for each $r$ such that $\beta_G(r)\ge\lambda$ we define a random graph $G_{r,\lambda}$ with vertex set $V$ by connecting vertices $u,v\in V$ by an edge with probability $\lambda/\beta_G(r)$ if $d_G(u,v)\le r$ and with probability $0$ otherwise, with the events that different edges are present independent. Equivalently, $G_{r,\lambda}$ is $G_r$ with a Bernoulli-$(\lambda/\beta_G(r))$ edge configuration. We will often call $r$ the \emph{scale} and $\lambda$ the \emph{degree parameter} of this model, although the expected degree of a vertex of $G_{r,\lambda}$ is actually $\frac{\beta_G(r)-1}{\beta_G(r)}\lambda$, not $\lambda$. In light of \eqref{eq:pc>}, \cref{thm:pc} is then easily seen to be equivalent to the following result.
\begin{theorem}\label{transitive}
Let $G$ be a transitive graph with superlinear polynomial growth. For any fixed degree parameter $\lambda > 1$, there exists $R>0$ such that if $r\ge R$ then $G_{r,\lambda}$ has an infinite component with probability $1$.
\end{theorem}
The proof of \cref{transitive} is broadly in four stages.

\subsection{Reduction to Cayley graphs of nilpotent groups}
We begin addressing \cref{transitive} by reducing the problem to the case of a Cayley graph of a nilpotent group. To achieve this we make use of a theorem of Trofimov \cite{trofimov} that describes the structure of an arbitrary transitive graph of polynomial growth and shows that such a graph $G$ can be approximated by a Cayley graph. The Cayley graph approximating $G$ also has polynomial growth. A celebrated theorem of Gromov \cite{Grm} then implies that the group underlying this Cayley graph is \emph{virtually nilpotent}, i.e. has a nilpotent subgroup of finite index. It is a standard fact (proven in \cref{virt.nilp->normal} below) that we may assume this nilpotent subgroup to be torsion-free. To conclude we obtain that there exists a quotient $G'$ of $G$ that admits a free action by a finitely generated torsion-free nilpotent group $\Gamma$, and our aim is to compare spread-out percolation on $G$ to spread-out percolation on a Cayley graph of $\Gamma$ defined with respect to a certain natural generating set.

Whilst this is a conceptually obvious thing to try, there are some further technical obstacles to overcome. One challenge is that, although a famous result of Benjamini and Schramm \cite[Theorem 1]{bperc96} allows us to compare Bernoulli-$p$ percolation on $G$ to Bernoulli-$p$ percolation on $G'$, this is of limited use in our context, since spread-out percolation with parameter $\lambda$ at scale $r$ corresponds to very different values of $p$ in Bernoulli percolation on $G'$ or on $G$. There is a similar issue in comparing spread-out percolation on $G'$ to spread-out percolation on $\Gamma$.

It turns out that we can adapt Benjamini and Schramm's argument to overcome this issue, but at the expense of comparing spread-out percolation at scale $r$ on $G$ to spread-out percolation at some other scale $r'$ on $G'$. We then need to compare spread-out percolation at different scales on $G'$, for which we use a quantitative version of a result of Tessera \cite[Theorem 4]{tessera.folner} comparing the sizes of balls of different radii in transitive graphs of polynomial growth.

As such, our proof of \cref{transitive} is reduced to the class of Cayley graphs, so it will be convenient to introduce some notation specialized to that setting.
If $\Gamma$ is a group with a finite symmetric generating set $S$, then we define the \emph{word length} $\Vert \gamma\Vert_S$ of an element $\gamma\in\Gamma$ to be the minimal integer $n\ge0$ such that $\gamma=s_1\cdots s_n$ for some $s_i\in S$. The \emph{word metric} $d_S$ on $\Gamma$ is then defined by $d_S(\gamma_1,\gamma_2)=\Vert \gamma_1^{-1}\gamma_2\Vert_S$. The \emph{Cayley graph} $G(\Gamma,S)$ has vertex set $\Gamma$ and edge set $\{(\gamma_1,\gamma_2)\in G\times G:\gamma_1 \in\gamma_2  S,\gamma_1\ne \gamma_2\}$. Note that the graph distance between two vertices $\gamma_1,\gamma_2$ of $G(\Gamma,S)$ is then exactly $d_S(\gamma_1,\gamma_2)$. Given $\gamma\in\Gamma$ and $r\ge0$, write $B_{\lVert \cdot \rVert_S}(\gamma,r)=\{\gamma'\in\Gamma:d_S(\gamma,\gamma')\le r\}$, the ball of radius $r$ centred at $\gamma$ in $G(\Gamma,S)$. Given a parameter $\lambda\ge0$, for each $r$ such that $\#B_{\lVert \cdot \rVert_S}(\id,r)\ge\lambda$ we write $G^{(S)}_{r,\lambda}(\Gamma)$ for the spread-out percolation on $G(\Gamma,S)$. Essentially $G^{(S)}_{r,\lambda}(\Gamma)$ is the graph with vertex set $\Gamma$ where we connect elements $\gamma_1,\gamma_2$ by an edge with probability $\lambda/\#B_{\lVert \cdot \rVert_S}(\id,r)$ if $d_S(\gamma_1,\gamma_2)\le r$ and with probability $0$ otherwise, with the events that different edges are present independent.

\cref{transitive} is then reduced to the following theorem.
\begin{theorem}\label{4.1}
Let $\Gamma$ be a torsion-free nilpotent group that is not isomorphic to $\mathbb{Z}$, and let $S$ be a finite symmetric generating set. For any fixed degree parameter $\lambda > 1$, there exists $R>0$ such that if $r\ge R$ then $G^{(S)}_{r,\lambda}(\Gamma)$ has an infinite component with probability 1.
\end{theorem}
It follows from Gromov's theorem (or an earlier result of Justin \cite{justin}) that the condition $\Gamma\not\cong\Z$ appearing in \cref{4.1} is equivalent in the class of torsion-free nilpotent groups to the condition that the growth is superlinear, so \cref{4.1} really is a special case of \cref{transitive}.

\subsection{Reduction to a graph in the asymptotic cone}
As mentioned in the introduction, we adapt an alternative proof of Penrose's result by Bollobás, Janson, and Riordan \cite{BJR}, which builds on their earlier work on inhomogeneous random graphs \cite{BJR2}.
Such an approach is not possible in an arbitrary graph, as it requires an underlying metric space with a measure.
The second stage of our proof establishes such a setting. Before we delve into the details, let us first provide a heuristic idea of the probabilistic argument contained in \cite{BJR}.

For simplicity, we consider the case where the graph $G$ is the square lattice $\mathbb{Z}^2$, equipped with a distance induced by a norm  $\Vert \cdot \Vert$ on $\mathbb{R}^2$. Set $\beta_r=\#\{x\in\Z^2:0<\|x\|\le r\}$. Given a fixed $\lambda>0$ and $r$ such that $\beta_r>\lambda$, the random graph $G_{r,\lambda}$ is obtained by connecting independently every two vertices at distance $r$ with probability $\lambda/ \beta_r$. Consider the graph $X_r$ with vertex set $(1/r)\cdot \mathbb{Z}^2$, and two vertices $x,y\in (1/r)\cdot \mathbb{Z}^2$ connected by an edge if and only if $\Vert x-y\Vert\leq 1$. The random subgraph $X_{r,\lambda}$ of $X_r$ is defined as $X_r$ with a Bernoulli-$(\lambda/\beta_r)$ edge configuration. The construction $X_{r,\lambda}$ has the same distribution as $G_{r,\lambda}$. This trivial restatement has the advantage of enabling us to fix a kernel in $\mathbb{R}^2$ that determines the percolation configuration for every $r$.

In \cite{BJR}, the authors, in order to study $X_{r,\lambda}$, begin by tessellating $\mathbb{R}^2$ into boxes $B = [a,b] \times [c,d]$ of uniform size. They define $\mathcal{S}$ as either a single box or the union of two neighboring boxes. The distribution of $ X_{r,\lambda}\cap \mathcal{S}$ is treated as a finite inhomogeneous random graph. Using results from \cite{BJR}, they conclude that for a sufficiently large box size, as $r$ grows to infinity, the following event holds with high probability: the largest connected component of $X_{r,\lambda}\cap \mathcal{S}$ has size $\Omega(\#(X_{r }\cap \mathcal{S}))$, while the second largest component has size $o(\#(X_{r }\cap \mathcal{S}))$. To derive the desired result from this local structure, the authors apply a renormalization argument. They construct a new graph whose vertex set consists of the boxes, two vertices are connected if and only if the corresponding boxes are neighboring, each contains connected components of linear size, and their union contains a unique connected component of linear size. This new random graph is a $2$-independent Bernoulli percolation model on $\mathbb{Z}^2$. The proof is completed by applying a classical result of Liggett, Schonmann and Stacey \cite{Lig} on $K$-independent percolation measures.

We will reduce the proof of \cref{4.1} to the case of another random graph where we can adapt the aforementioned argument. For this, we first require an analogue of the embedding of $\mathbb{Z}^d$ into $\mathbb{R}^d$ for a torsion-free nilpotent group $\Gamma$. We will consider $\Gamma$ as a discrete subset of its \emph{asymptotic cone} $G_\infty$, a notion defined by Gromov \cite{Grm}, a Lie group that is a \emph{Gromov--Hausdorff limit} (see \cite{Grm2}) of the sequence of metric spaces $(\Gamma,d_S/n)$. Pansu \cite{Pan} described $G_\infty$ and equipped it with a metric $d_{CC}$, now called the \emph{Pansu limit metric}, which asymptotically behaves like the word metric for the elements of the group $\Gamma$. Additionally \cite{Pan} proves the existence of  
a constant $c_S>0$ -- the \emph{coefficient of polynomial growth} -- such that $\#B_{\lVert \cdot \rVert_S}(\id,n)= c_S n^{d_\Gamma}+ o(n^{d_\Gamma})$ for all $n\in\N$, where $d_\Gamma$ is a non-negative integer -- the \emph{degree of polynomial growth} -- given by the Bass-Guivarc'h formula \cite{Gui}. (See \cref{4.1.2} below for a more precise estimate of the error term due to Breuillard and Le Donne).

For $\lambda\geq 0$ and $r$ such that $(c_S r^{d_\Gamma}) \geq \lambda$, we define the new random graph $G^{r,\lambda}=G^{r,\lambda}(\Gamma,S,\mathrm{d}_{CC})$ with vertex set $\Gamma$. Two vertices $\gamma_1,\gamma_2$ are connected by an edge with probability $ \lambda/(c_S r^{d_\Gamma})$ if $d_{CC}(id, \gamma_2^{-1} \gamma_1)\leq r$, and with probability $0$ otherwise. Essentially, the graph $G^{r,\lambda}$ is a model of spread-out percolation on the elements of $\Gamma$, but defined with respect to the Pansu limit metric.

\begin{theorem}\label{4.2}
Let $\Gamma$ be a torsion-free nilpotent group that is not isomorphic to $\mathbb{Z}$, and let $S$ be a finite symmetric generating set. 
For any fixed degree parameter $\lambda >1$, there exists $R>0$ such that if $r\geq R$, then $G^{r,\lambda}(\Gamma,S,\mathrm{d}_{CC})$ has an infinite component with probability 1.
\end{theorem}

Breuillard and Le Donne \cite{BrL} have provided sharp bounds for the difference between the word metric and the Pansu limit metric (see the precise statement in \cref{4.1.1}). These bounds make the comparison of the two random graphs $G^{r,\lambda}=G^{r,\lambda}(\Gamma,S,\mathrm{d}_{CC})$ and $G^{(S)}_{r,\lambda}(\Gamma)$ straightforward. The advantage of such a comparison is that it allows us to reduce the proof of \cref{4.1} to the case of the random graph $G^{r,\lambda}$, which is defined within a metric space equipped with a measure.

However, the product $\gamma_2^{-1}\gamma_1$ in the definition of $G^{r,\lambda}$ corresponds to the group operation of $\Gamma$, which does not necessarily coincide with the group operation of $G_\infty$: the group $\Gamma$ is not generally a sub\textit{group} of the asymptotic cone, but merely a sub\textit{set}. Consequently, the argument of Bollobás, Janson and Riordan cannot be applied, which leads us to the next stage of the proof.

\subsection{Showing $G^{r,\lambda}$ is locally an inhomogeneous random graph} Our goal is to study $G^{r,\lambda}$ locally, treating it as an inhomogeneous random graph. It is not feasible to construct a fixed kernel that determines the percolation configuration for $G^{r,\lambda}$ as $r$ varies. This limitation arises because the dilations do not preserve the structure of $\Gamma$. We address this issue by considering a graphical sequence of kernels (definition in \cref{4}) that converge to a limit kernel defined in $G_\infty$. To achieve this, we require a Lie group $G$ that respects the group structure of $\Gamma$. This Lie group $G$ is given by a classical result of Malcev \cite{Mal}, who proved that $\Gamma$ embeds as a lattice in a nilpotent Lie group $G$, called the \emph{Malcev closure} of $\Gamma$. However, the Malcev closure generally does not admit any natural rescaling function. To address these problems, we adapt an idea of Breuillard and Le Donne \cite{BrL}, who treat $G$ and $G_\infty$ as one manifold with two different group structures, with operations $\ast$ and $\ast_\infty$ respectively, and consider $\Gamma$ as a discrete set in this manifold. We construct a sequence of kernels in $G$ with respect to the $\ast$-operation, which determines the percolation configuration for $G^{r,\lambda}$, and converges to a kernel described in terms of the $\ast_\infty$-operation.

A key ingredient for this approach is a property of the Lebesgue measure: for any measurable subset $A\subset\R^d$ with $\mu(\partial A)=0$, the ratio $\#((1/r)\cdot \Z^d\cap A)/r^d$ converges to the Lebesgue measure $\mu(A)$ (this statement is equivalent to the one presented in the introduction, see \cite{Patrick}). We establish an analogue of this result, which is a purely deterministic theorem concerning the Haar measure on the asymptotic cone $G_\infty$ of the torsion-tree nilpotent group $\Gamma$.

For any $\rho>0$ the asymptotic cone $G_\infty$ can be endowed with a rescaling function $\delta_\rho$ referred to as \textit{dilation}. Moreover, $G_\infty$ is a Lie group with a Haar measure $\mathcal{H}$, normalized so that $\mathcal{H}$ assigns volume $1$ to the unit ball with respect to the Pansu limit metric. A subset $A\subset G$ is called an \emph{$\mathcal{H}$-continuity set} if it is $\mathcal{H}$-measurable and satisfies $\mathcal{H}(\partial A)=0$.
\begin{theorem}\label{Haartheo}
Let $\Gamma$ be a torsion-free nilpotent group with a finite symmetric generating set $S$, and let $A\subset G_\infty$ be a bounded $\mathcal{H}$-continuity set. Then,
$$
\lim_{r\to \infty}\frac{\# \left( \delta_{1/r}(\Gamma)\cap A\right)}{c_S r^{d_\Gamma}}=\mathcal{H}(A),$$
where $c_S$ is the coefficient of polynomial growth.
\end{theorem}
Note that the Haar measure $\mathcal{H}$ depends on the generating set $S$, as its normalization to assign volume $1$ to the unit ball relies on the choice of generating set. 
Theorem \ref{Haartheo} can be interpreted measure-theoretically as follows. Define a discrete measure $\nu_r$ on $G_\infty$, assigning mass $1/c_S r^{d_\Gamma}$ to each element of $\delta_{1/r}(\Gamma)$. Let $A$ be a bounded $\mathcal{H}$-continuity subset of $G$. Then, $\nu_r(A)\to\mathcal{H}(A)$ as $r\to\infty$. This implies that if we restrict to a bounded subset, $\nu_r$ converges weakly to $\mathcal{H}$.

\subsection{Renormalization argument}
We finally consider the intersection of the graph $G^{r,\lambda}$ with some suitable region $\mathcal{S}$ of $G$, whose role is analogous to that of the square boxes from the Bollobás--Janson--Riordan argument. 
However, while in $\R^d$ it is possible to tessellate the space with isometric copies of a box indexed by the elements of $\Z^d$, there is no analogous tessellation of $G$ by isometric copies of $\mathcal{S}$ indexed by elements of $\Gamma$. One key reason for this is that the metric $d_{CC}$ is not $\ast$-invariant.
Moreover, the region $\mathcal{S}$ must be defined in a way such that it can grow arbitrary large, depending on how close to $1$ the degree parameter $\lambda$ is.
To address this, for every size of the region $\mathcal{S}$, we find a subgroup $H$ of $\Gamma$ isomorphic to $\Z^2$ that is also a subgroup of $G_\infty$, i.e., a subgroup of $G$ on which the $\ast$-operation and the $\ast_\infty$-operation coincide. However, this still leaves another issue: the translates of $\mathcal{S}$ by elements of $H$ do not fit together as nicely as the boxes in $\R^d$, because the different translates overlap one another. As the size of the region grows, the number of overlaps changes. To ensure that this does not ruin the argument, we must define the region $\mathcal{S}$ and the subgroup $H$ carefully so that the number of these overlaps is uniformly bounded. This will allow us to define a $K$-dependent percolation on $H$ for some fixed $K$.

\begin{remark}
Since this work is within the same framework as that of Contreras, Martineau and Tassion \cite{CMT.supercritical} and \cite{CMT.locality}, we should make the following remark. The final step in the proof of \cref{4.2} is a renormalization argument, where we consider a Bernoulli percolation model on the square lattice. Each vertex $v$ of $H\cong\Z^2$  corresponds to the translation of $\mathcal{S}$ by $v$. An edge is considered open if and only if each corresponding region contains a connected component of linear size and their union contains exactly one connected component of linear size. The proof is then completed by a classical result of Liggett, Schonmann, and Stacey \cite{Lig} on $K$-independent percolation measures. It might be possible to replace this final step by using the tools developed in the work of Contreras, Martineau, and Tassion \cite{CMT.locality}, which also contains a renormalization argument and the result from \cite{Lig}. In that case, the treatment would focus on the balls of radius $n$ for the spread-out percolation on $\Gamma$. However, this would not shorten the proof, as in order to establish the asymptotic properties, we must fix a region in $G_\infty$, with the ball of radius $n$ with respect to $\mathrm{d}_{CC}$ being the most natural choice. Working with the balls of $G_\infty$ is just as complicated as working with the regions $\mathcal{S}$ that we considered. Thus, adopting such an approach would not significantly alter the proof. 
\end{remark}

\begin{remark}
The results of Penrose and Bollobás, Janson and Riordan are actually rather more general than what we have described, in that instead of considering the norm $\|\,\cdot\,\|$ one can consider an arbitrary symmetric probability distribution $\phi$ on $\R^d$, and put an edge between $x,y\in\Z^d$ with probability
\[
\frac{\lambda\phi(x-y)}{\sum_{z\in\frac1r\Z^d\setminus\{0\}}\phi(z)}.
\]
The version we stated above is of course the special case of this in which $\phi$ is the characteristic function of the closed unit ball of the norm.

In light of the issues just discussed concerning the embedding of $\Gamma$ into its Malcev closure, there does not appear to be a natural analogue of this generalization in the setting of \cref{4.1}. There certainly does not appear to be any such natural analogue in the setting of \cref{thm:pc}.
\end{remark}

\section{Preliminaries on nilpotent groups}\label{2}
\subsection{Lattices in nilpotent Lie groups}
Let $\Gamma$ be a group. If $\Gamma^{(s)}\neq \Gamma^{(s+1)}=1$, then $\Gamma$ is called $s$-step nilpotent.
If moreover $\Gamma$ is finitely generated, the quotients $\Gamma^{(i)}/\Gamma^{(i+1)} $ are Abelian groups, we denote $d_i$ their rank. The constant $d_\Gamma$ given by the Bass-Guivarc'h formula $d_\Gamma = \sum_i i d_i$ is the growth degree of the group (see \cite{Gui}), i.e. if $S$ is a finite, symmetric set of generators containing the identity, then there are $a,b,c>0$ such that $ar^{d_\Gamma}< \# \vert S^r \vert<br^{d_\Gamma}$ for every $r>c$. 

Let $\Gamma$ be a finitely generated, torsion-free, nilpotent group, there exists a connected, simply connected, nilpotent Lie group $G$ such that $\Gamma$ is isomorphic to a lattice in $G$, we call $G$ the Malcev closure of $\Gamma$ (see \cite{Rag, Mal}).
A discrete subgroup $\Gamma<G$ is called a lattice in $G$ if $G/\Gamma$ carries a finite invariant measure, in this case that $G$ is nilpotent, $\Gamma$ being a lattice is equivalent to $G/\Gamma$ being compact (see \cite{Rag}). From now on all nilpotent Lie groups will be assumed connected and simply connected.

Let $G$ be a nilpotent Lie group.
We denote by $\mathfrak{g}$ the Lie algebra of $G$. The descending central series of $\mathfrak{g}$ is defined as $\mathfrak{g}^{(i+1)}=[\mathfrak{g},\mathfrak{g}^{(i)}]=\mathbb{R}\text{-span}\lbrace[X,Y]:X\in \mathfrak{g}, Y\in \mathfrak{g}^{(i)}\rbrace $, since $\mathfrak{g}^{(i)}$ is the Lie algebra of $G^{(i)}$ for every $i=1,\ldots, s$ (see \cite{Hoc}), then $\mathfrak{g}$ is a nilpotent Lie algebra.
There exists a basis $\lbrace X_1,\ldots,X_d\rbrace$ for $\mathfrak{g}$ such that $\mathfrak{h}_n=\mathbb{R}\text{-span}\lbrace X_1,\ldots,X_n\rbrace$ is an ideal of $\mathfrak{g}$ and $\mathfrak{h}_{n_i} = \mathfrak{g}^{(i)}$,
for every $i=1,\ldots,s$, this base is a strong Malcev basis for the descending central series (see \cite{Grn}). 

The exponential map $\exp: \mathfrak{g}\to G$ is defined globally and it is an analytic diffeomorphism, moreover the Baker-Campbell-Hausdorff formula (BCH formula, in abbreviation) $$\exp X \ast \exp Y = \exp( X+Y+\frac{1}{2}[X,Y]+\frac{1}{12}[X,[X,Y]]-\frac{1}{12}[Y,[X,Y]]\ldots)\label{BHC}$$
holds for every $X,Y\in\mathfrak{g}$ (see \cite{Grn}), by $\ast$ we denote the group product in $G$. We denote $\log:=\exp^{-1}$ and $X\bullet Y:=\log(\exp X \ast \exp Y)$ for every $X,Y\in \mathfrak{g}$.

There are two coordinate systems that Malcev endowed $G$ with, the first one is the exponential coordinates, $\Phi:\mathbb{R}^d\to G$, defined as $\Phi(x_1,\ldots,x_d)=\exp(x_1X_1+\ldots + x_d X_d)$ for $(x_1,\ldots,x_d)\in \mathbb{R}^d$, also called canonical coordinates of the first kind. The second one is $\Psi:\mathbb{R}^d\to G$, defined as $\Psi(x_1,\ldots,x_d)=\exp(x_1X_1)\ast\ldots\ast\exp(x_dX_d)$ called strong Malcev coordinates or canonical coordinates of second kind.

We will identify $\mathfrak{g}$ with $\mathbb{R}^d$ through the linear isomorphism $(x_1,\ldots,x_d)\mapsto x_1X_1+\ldots x_dX_d$, we will abuse notation and write elements of $\mathfrak{g}$ as vectors $X=(x_1,\ldots,x_d)$, then the exponential coordinates will be written as $\Phi(x_1,\ldots,x_d)=\exp X$. There is always a Malcev basis for the descending central series such that $\Gamma = \lbrace \Psi(\vec{v}) :\vec{v} \in\mathbb{Z}^d\rbrace$ (see \cite{Matsu}), this base is called strong Malcev basis for the descending central series strongly based on $\Gamma$.

Both $\Phi$ and $\Psi$ push forward the Lebesgue measure from $\mathbb{R}^d$ to a Haar measure on $G$ (see \cite{Grn})\label{LebHar}, moreover the following relation holds; if $\mu$ is the pushforward from $\Phi$ and $\nu$ is the pushforward from $\Psi$ then

$$\nu=\vert \det \mathbb{J}_{(\Psi^{-1}\circ \Phi)}(0)\vert \mu\label{DetHaar}.$$

\subsection{Carnot groups}
We will describe the structure of a Carnot group, a lot of the material of this section is based on a survey on Carnot groups \cite{LeDonne1}. 
Let $\mathfrak{g}$ be a Lie algebra, we say that $\mathfrak{g}$ admits a stratification, if as a vector space it can be decomposed as $\mathfrak{g}=V_1\oplus \ldots \oplus V_s$, for $s\in \mathbb{N}$, such that $V_s\neq \lbrace 0 \rbrace, [V_1,V_j]=V_{j+1}$, for $1\leq j \leq s-1$ and $[V_1,V_s]=\lbrace 0 \rbrace$. The vector subspaces $V_j\subset \mathfrak{g}$ are called strata of the stratification. If there is associated a fixed stratification for $\mathfrak{g}$ then the Lie algebra is called stratified. By definition a stratified Lie algebra is nilpotent of nilpotency step $s$. The first stratum $V_1$ determines the stratification and different stratifications are linearly isomorphic.
Let $G$ be a connected simply connected Lie group, we say that $G$ is stratified (admits stratification) if the associated Lie algebra $\mathfrak{g}$ is stratified (admits stratification correspondingly).

A Carnot group is a stratified Lie group enhanced with a subFinsler structure with distribution induced by the first stratum of the  associated stratification. We will provide the definition of a subFinsler space, but in order to keep the notions as simple as possible, we will continue in the setting of a connected, simply connected Lie group. Though the general definition of a subFinsler manifold does not require for the manifold to have a Lie group structure.

Let $V\subset \mathfrak{g}$ be a linear subspace and fix $\lVert \cdot \rVert$ a norm on it, then $V$ determines a distribution $\Delta$, meaning a subbundle, of the tangent bundle of $G$. For every $g\in G$ we obtain $\Delta_g := (L_g)^\ast V$, where $L_g$ is the left translation on $G$ and $F^\ast$ is the pull back operator of a diffeomorphism $F$. The norm on $V$ induces a norm on $\Delta_g$, $\Vert X_g \Vert = \Vert (L_g)^\ast X_e \Vert:= \Vert X_e \Vert$, for every $X_g\in \Delta_g$ and $g\in G$. If the distribution $\Delta$ is bracket generating, meaning that the Lie algebra generated, with respect to the Lie bracket, by the vector fields $\Gamma(\Delta)$ is the whole algebra $\mathfrak{g}$, then we call $(G,\Delta,\lVert \cdot \rVert)$ subFinsler manifold. If the norm $\lVert \cdot \rVert$ is induced from an inner product, then it is called subRiemannian manifold.

The bracket generating distribution $\Delta$ is called horizontal and an absolutely continuous curve $\gamma:[a,b]\to G$ is also called horizontal or admissible if its derivative $\dot{\gamma}(t)\in \Delta_{\gamma(t)}$ for almost all $t\in [a,b]$. Every such a curve can be associated with a length $\text{Length}_{\lVert \cdot \rVert}(\gamma)=\int_a^b\Vert \dot{\gamma}(t)\Vert \mathrm{d} t$, since $\Vert \dot{\gamma}(t)\Vert$ is defined almost for every $t$. A subFinsler manifold is endowed with a distance function, for every two points $g_1,g_2 \in G$ the distance is defined as 
$d_{CC}(g_1,g_2)=\inf \lbrace \text{ Length }_{\lVert \cdot \rVert}(\gamma):\gamma \text{ admissible curve from } g_1 \text{ to } g_2\rbrace$. It is
called Carnot-Carathéodory distance or subFinlser distance.

From Chow's Theorem (see \cite{Mon}) the Carnot-Carathéodory metric is finite, induces the manifold topology and it is geodesic, i.e. that for every $g_1,g_2\in G$ there is an admissible $\gamma$ from $g_1$ to $g_2$ such that $d_{CC}(\gamma(t_1),\gamma(t_2))=\vert t_1 - t_2 \vert$ for every $t_1,t_2\in [a,b]$. If $\gamma$ is an admissible curve, then $g\dot{\gamma}(t)\in\Delta_{g\gamma(t)}$ and $\text{Length}_{\lVert \cdot \rVert}(\gamma)=\text{Length}_{\lVert \cdot \rVert}(g\gamma)$ since $\Vert \dot{\gamma}\Vert = \Vert \dot{g\gamma}(t)\Vert$ by definition, $d_{CC}(gg_1,gg_2)=d_{CC}(g_1,g_2)$. 

Let $G$ be a Carnot group with Lie algebra $\mathfrak{g}$, it can be proved that $\mathfrak{g}^{(k)}=V_k \oplus \ldots \oplus V_s$ for every $k=1,\ldots,s$, where $V_i$ are the strata of the associated stratification. Indeed one needs to observe that $[V_i,V_j]\subset V_{i+j}$ for every $i,j=1,\ldots,s$, setting $V_k=\lbrace 0 \rbrace$ for $k>s$.

\begin{note}
Every nilpotent Lie algebra could be associated with a graded Lie algebra but the linear isomorphism defined above does not necessarily extend to a Lie algebra isomorphism, it does extend only when the group is stratified.
\end{note}

\paragraph{Dilations}
For every $\lambda\in \mathbb{R}$ we set $\delta_{\lambda}^{\mathfrak{g}}=\lambda^i X$ for every $X\in V_i$ and $i=1,\ldots,s$. This map extends linearly to a map $\delta_\lambda^{\mathfrak{g}}:\mathfrak{g}\to \mathfrak{g}$ called algebra dilation. If $\lambda \neq 0$, then $\delta_\lambda$ is an automorphism of $\mathfrak{g}$. Since $G$ is connected and simply connected, every $\phi$ Lie algebra automorphism induces a unique $F$ Lie group automorphism, such that $dF=\phi$ and $\exp\circ dF = F\circ \exp$. For every $\lambda\neq 0$ the induced from the dilations automorphisms $\delta^G_{\lambda}:G\to G$ are called intrinsic dilations and they satisfy $\exp \circ \delta_\lambda^{\mathfrak{g}}=\delta_\lambda^G\circ \exp$. Dilations could be equivalently defined as $\delta_\lambda^G(g):=\exp \circ \delta_\lambda^{\mathfrak{g}}\circ \exp^{-1}(g)$ for every $g\in G$. Some authors use the notation $\log:=\exp^{-1}:G\to \mathfrak{g}$, moreover for simplicity we will write $\delta_\lambda$ in both cases.

Let $G$ be a Carnot group and $\delta_\lambda$ the group dilations, then for any $X\in V_1$ we have defined $\delta_\lambda(X)=\lambda X\in V_1$, this means that $\Vert \delta_\lambda(X)\Vert = \lambda\Vert X\Vert$. Hence for any admissible curve $\gamma$ that connects $g_1$ with $g_2$ in $G$, we obtain that $\delta_\lambda\circ \gamma$ connects $\delta_\lambda(g_1)$ with $\delta_\lambda(g_2)$ and $\frac{d}{dt}(\delta^{G}_\lambda\circ \gamma)(t)=
d\delta^G_\lambda\circ \dot{\gamma}(t)=\delta_\lambda^{\mathfrak{g}}(\dot{\gamma}(t))=\lambda\dot{\gamma}(t)\in V_1$. 
Therefore $\delta_\lambda\circ\gamma$ is admissible and $\text{Length}_{\lVert \cdot \rVert}(\delta_\lambda\circ \gamma)=\lambda \text{Length}_{\lVert \cdot \rVert}(\gamma)$, this implies that $d_{CC}(\delta_\lambda(g_1),\delta_\lambda(g_2))=\lambda d_{CC}(g_1,g_2)$ and $B_{CC}(g,r)=\delta_r(B_{CC}(g,1))$. Lastly from the Baker-Campbell-Hausdorff formula (\ref{BHC}) and the linearity of $\delta_\lambda^{\mathfrak{g}}$, we obtain that $\delta_\lambda(g_1)\delta_\lambda(g_2)=\delta_\lambda(g_1g_2)$, for every $g_1,g_2\in G$.
We will denote $B_\infty:=B_{CC}(\id,1)$.

\subsection{Graded algebra}
We will provide a description of the graded algebra, we adapt material from \cite{BrL} in our setting.
If $G$ is a connected, simply connected nilpotent Lie group, even when it is not stratifiable, there is a stratified Lie algebra arising from its Lie algebra. Let $\mathfrak{g}$ be the Lie algebra of $G$ and $\mathfrak{g}^{(n)}=[\mathfrak{g},\mathfrak{g}^{(n-1)}]$ its descending series. The direct sum $\mathfrak{g}_\infty :=
\oplus_{i=1}^s \mathfrak{g}^{(i)}/\mathfrak{g}^{(i+1)}$ endowed with a Lie bracket $[\cdot,\cdot]_{\infty}:\mathfrak{g}_\infty\times\mathfrak{g}_\infty\to \mathfrak{g}_\infty$ obtained from the following relations; let $[X]_i\in \mathfrak{g}^{(i)}/\mathfrak{g}^{(i+1)}$ and $[Y]_j\in \mathfrak{g}^{(j)}/\mathfrak{g}^{(j+1)}$, then $[[X]_i,[Y]_j]_\infty=[[X,Y]]_{i+j}\in \mathfrak{g}^{(i+j)}/\mathfrak{g}^{(i+j+1)}$, is a stratified Lie algebra called the graded algebra of $\mathfrak{g}$. The relations are well defined and determine uniquely a stratified Lie structure on $\mathfrak{g}_\infty$. If the Lie algebra $\mathfrak{g}$ is not stratifiable then obviously $\mathfrak{g}$ and $\mathfrak{g}_{\infty}$ are different Lie structures. The stratified Lie algebra $\mathfrak{g}_\infty$ has an associated connected, simply connected nilpotent Lie group $G_\infty$ called the graded group of $\mathfrak{g}$.

Let $\lbrace X_1,\ldots, X_d\rbrace$ be a strong Malcev basis for the descending series $\mathfrak{g}\supset \mathfrak{g}^{(2)}\supset  \ldots \supset  \mathfrak{g}^{(s)}$. For every $j=1,\ldots s$ 
$\mathfrak{g}^{(j)}=\mathfrak{g}^{(j+1)}\oplus \left<X_{m_{s-j}+1},\ldots, X_{m_{s-j+1}}\right>$. Set $V_j:= \left<X_{m_{s-j}+1},\ldots, X_{m_{s-j+1}}\right>$, as a vector space, the projection $p_{j+1}: \mathfrak{g}^{(j)}=\mathfrak{g}^{(j)}/\mathfrak{g}^{(j+1)}$ is restricted to a linear isomorphism $L_j:V_j\to\mathfrak{g}^{(j)}/\mathfrak{g}^{(j+1)}$, which consecutively extends to $L$ a linear isomorphism $L:\mathfrak{g}\to \mathfrak{g}_{\infty}$. 

The exponential map $\exp : \mathfrak{g}_{\infty} \to G_\infty$ is a diffeomorphism and so it is the linear map $L^{-1}:\mathfrak{g}_{\infty} \to \mathfrak{g}$, thus the composition 
$G_\infty\xrightarrow{\log} \mathfrak{g}_{\infty}\xrightarrow{L^{-1}}\mathfrak{g}\xrightarrow{\exp} G$ is a diffeomorphism. This map is not necessarily group homomorphism 
and if the group $G$ is not stratifiable then it is certainly not. Through this map though we can endow $G$ with a new product rule from $G_\infty$, we denote it by $\ast_\infty$, $\mathfrak{g}$ correspondingly receive a new Lie bracket $\left[\cdot,\cdot\right]_\infty$ and a new product $\bullet_\infty$. We denote the induced canonical coordinates of second kind by $\Psi_\infty:\mathbb{R}^d \to G$.

Let $\Gamma$ be a discrete, torsion-free, nilpotent group and $G$ its Malcev closure. When we follow this procedure for the stratification of $G$, we choose $\lbrace X_1,\ldots, X_d\rbrace$ to be a strong Malcev basis that is strongly based on $\Gamma$. To clarify things $\Gamma$ is a subgroup of $(G,\ast)$ but $\Gamma$ is not in general a subgroup of $(G,\ast_\infty)$; $\Psi(\mathbb{Z}^d)=\Gamma$ but $\Psi_\infty(\mathbb{Z}^d)$ is not even necessarily a group. 

\paragraph{Pansu limit metric.}
We will present a description of the Pansu limit metric \cite{Pan}, thus of the asymptotic cone. 
Let $S$ be a symmetric set of generators for $\Gamma$.
From the BCH formula (\ref{BHC}), the following map $\pi : G\xrightarrow{\log} \mathfrak{g}\xrightarrow{pr} \mathfrak{g}/[\mathfrak{g},\mathfrak{g}]\to V_1$ is an homomorphism, $\pi (g_1\ast g_2)=\pi(g_1)+\pi(g_2)$.
Then $\pi(\Gamma)$ is a lattice in $V_1$ and $\pi(S)$ generates it. We can endow $\mathfrak{g}/[\mathfrak{g},\mathfrak{g}]$ with a norm $\lVert \cdot\rVert_Q$ such that the unit ball $Q$ of the norm is the convex hull of $\pi (S)$.

Since $G$ endowed with $\ast_\infty$ is a stratified Lie group with the first stratum of its Lie algebra being $\mathfrak{g}/[\mathfrak{g},\mathfrak{g}]$, the triple $(G,V_1,\lVert \cdot \rVert_{Q})$ induces $d_{CC}$, a Carnot-Carathéodory metric on $G$, called the Pansu limit metric and $G_\infty$ is a Carnot group.
We just mention that $G_\infty$ endowed with this metric is the asymptotic cone of $\Gamma$ with respect to the word metric arising from $S$. This metric is well comparable to the word metric in the sense that \cite{Pan} for every $\gamma\in\Gamma$ that $\Vert \gamma\Vert_S\to \infty$ the limit $\frac{\Vert \gamma \Vert_S}{d_{CC}(e,\gamma)}\to 1$. 

Since the graded algebra $\mathfrak{g}_\infty$ is stratified, we can define dilations on it and through the diffeomorphism pull back the dilations on $G$. For every $\lambda>0$ the dilations $\delta_{\lambda}^{\mathfrak{g}}:\mathfrak{g}\to \mathfrak{g}$ satisfy the relations $\delta_{\lambda}(X)=\lambda^j X$ for $X\in V_j$, specifically $\delta_{\lambda}^{\mathfrak{g}}(X_{i_j})=\lambda^j X_{i_j}$ for every $i_j=m_{{s-j}+1},\ldots,m_{s-j+1}$, then $\delta^{G}_{\lambda}:=\exp \circ \delta_{\lambda}\circ \log$.
Dilations do not behave well with the original group structure or the original Lie bracket, but this is natural because they are pulled back from $G_\infty$.

From a refinement of Pansu's construction, Breuillard and Le Donne \cite{BrL} proved that:
\begin{proposition}\label{4.1.1}
Let $\Gamma$ be a torsion-free, finitely generated, nilpotent group and $S$ a finite symmetric set of generators for $\Gamma$. If $\Vert \gamma \Vert_S\to \infty$, then
\begin{align}
\vert d_{CC}(\id,\gamma)-\Vert\gamma\Vert_S \vert = O\left(d_{CC}(\id,\gamma)^{\alpha_s}
\right)
\end{align}
where $\alpha_s=1-\frac{2}{3s}$ and $s$ is the nilpotency step.
\end{proposition}
They used this bound to provide estimates for the size of the ball:
\begin{proposition}\label{4.1.2}
Let $\Gamma$ be a torsion-free, finitely generated, nilpotent group and $S$ a finite symmetric  set of generators for $\Gamma$. Then
\begin{align}
\#\left(B_{\lVert \cdot \rVert_S}(\id,n)\right) =c_S n^{d_\Gamma}+O(n^{\beta_s})
\end{align}
where $\beta_s=d_\Gamma - \frac{2}{3s}$ and $c_S$ is a constant that depends on the set of generators.
\end{proposition}

\section{Convergence to Haar measure}\label{3}
In this section we prove Theorem \ref{Haartheo}. 
We begin by defining a non-uniform rescaling process in $\mathbb{R}^d$ that it will be pushed-forward to the dilations through the Malcev coordinates.

Let $\mu$ be the Lebesgue measure on $\mathbb{R}^d$. For every $\lambda >0$ set $\delta_{\lambda}:\mathbb{R}^d \to \mathbb{R}^d$ be the function
$(\vec{x_1},\vec{x_2},\ldots,\vec{x_n})\mapsto(\lambda\vec{x_1},\lambda^2\vec{x_2},\ldots,\lambda^n\vec{x_n})$, we view $\mathbb{R}^d$ as $V_1\oplus \ldots \oplus V_n$ and write $x\in \mathbb{R}^d$ as $(\vec{x_1},\ldots,\vec{x_n})$ with $\vec{x_i}\in V_i$. We define $d_i=\dim V_i$ and $d'=\sum_{i=1}^n id_i$, then $d=\sum_{i=1}^{n} d_i$. Recall that $A$ is a $\mu$-continuity set if the boundary has measure $\mu(\partial A)=0$ and the boundary $\partial A$ is defined as the elements in the closure of $A$ that do not belong in the interior of $A$.

\begin{lemma}\label{3.1}
For every $A\subset \mathbb{R}^d$ bounded $\mu$-continuity set,
$$\lim_{r\to\infty}\frac{ \# \left(\delta_{1/r} (\mathbb{Z}^d) \cap A\right)}{r^{d'}} 
=\mu (A).$$
\end{lemma}
\begin{proof}
Let $\epsilon >0$, we can find a $U_1$ open set such that $\mu (U_1)\leq \epsilon$ and $\partial A\subset U_1$, since $\mu(\partial A)=0$. Also from completeness of the Lebesgue measure there exists a compact set $K\subset \text{int} (A)$ with $\mu (A) \leq \mu(K)+\epsilon$. There are $U_2$ and $U_3$ such that $K\subsetneq U_2$, $(\text{int}(A))^c\subsetneq U_3$ and $U_2\cap U_3 =\emptyset$, this is from the $T_4$ property.

Let $x\in\mathbb{Z}^d$, write $x=(x^{(1)}_1,\ldots ,x^{(n-1)}_{d_{n-1}},x^{(n)}_{1},\ldots,x^{(n)}_{d_n})$, where $\vec{x_i}=(x^{(i)}_1,\ldots, x^{(i)}_{d_i})$ and set 
$C_{x}^k:= \left[x^{(1)}_1,x^{(1)}_1+2^{-k}\right]\times\ldots\times
\left[x^{(i)}_j,x^{(i)}_j+2^{-ik}\right]\times\ldots \times \left[x^{(n)}_{d_n},x^{(n)}_{d_n}+2^{-nk}\right]$ and \linebreak
$D_{x}^k:= \left(x^{(1)}_1-2^{-k},x^{(1)}_1+2^{-k}\right)\times\ldots\times  \left(x^{(n)}_{d_n}-2^{-nk},x^{(n)}_{d_n}+2^{-nk}\right)$ 
For every element $u\in U_2$ there exists $B_E(u,r)\subset U_2$ (the Euclidean ball in $\mathbb{R}^d$). From the density of diadic numbers we can find a $k>0$ and $x\in\delta_{2^{-k}}(\mathbb{Z})$ such that $u\in D_{x}^{k}\subset B_E(u,r)$.

We cover $U_2$ with $D_{x}^k\subset U_2$ then $K$ is also covered, which is compact, therefore it has a finite cover $\cup D_{x}^{k}$, without loss of generality we can assume that $k$ is the same for all $x$. 
Every $D_{x}^k \subset U_3^{c}$ since $D_{x}^k\subset U_2$ and therefore $\overline{D_{x}^k }\subset U_3^{c}$. We can write $\overline{D_{x}^k }$ as a union of $C^k_{y}$ that have foreign interiors.
Note that $x\in \delta_{2^{-k}}(\mathbb{Z})\cap A$ and $\mu (C^k_{x})=2^{-d'k}$.
Therefore we obtain the inequality 
$$\mu (A) - \epsilon \leq \mu (K) \leq \mu (\cup C^{k}_{x}) \leq 
\frac{ \#\left(\delta_{2^{-k}}(\mathbb{Z})\cap A\right)}{2^{-d'k}}.$$

We want to make this inequality valid for all $r>0$, not only for powers of $2$. Set the cardinality of the cover as $N_1 :=\#\left(\delta_{2^{-k}}(\mathbb{Z})\cap 
(\cup C^{k}_{x})\right)$, for every $r\in\mathbb{R}_{>0}$ and $x\in \delta_{r^{-1}}(\mathbb{Z})$ we define the sets 
$E_{x}^r:= \left[ x_1^{(1)},x_1^{(1)}+r^{-1}\right)\times\ldots\times
\left[ x_{d_n}^{(n)},x_{d_n}^{(n)} + r^{-n} \right)$
We can focus on the case of $r>2^k$. Every $C_{x}^k$ contains at least
$\Bigl\lfloor \frac{r}{2^k} \Bigl\rfloor^{d_1} \Bigl\lfloor \left( \frac{r}{2^k}\right)^2 \Bigl\rfloor^{d_2} \ldots \Bigl\lfloor \left(\frac{r}{2^k}\right)^n\Bigl\rfloor^{d_n}$ sets $E^r_{y}$, for $y\in \delta_{r^{-1}}(\mathbb{Z})$ and at most 
$\left(\Bigl\lfloor \frac{r}{2^k} \Bigl\rfloor+2\right)^{d_1} 
\left( \Bigl\lfloor  \frac{r}{2^k}\Bigl\rfloor^2 +2\right)^{d_2} \ldots \left( \Bigl\lfloor \frac{r}{2^k}\Bigl\rfloor^n+2\right)^{d_n}$.

Since $E_{y}^r$ have foreign interiors, if we consider the union of $E_{y}^r$ that are contained in $C_{x}^k$ then
$\mu(C_{x}^k)-\mu(\cup E_{y}^r)=\mu ( C_{x}^k\setminus 
(\cup E_{y}^r))< g(r)$, where
\begin{align*}
g(r):=&\sum_{i=1}^n d_i (2^{-k}+2r^{-1})^{d_1}\ldots
(2^{-nk}+2r^{-n})^{d_n}\frac{2r^{-d_i}}{(2^{-ik}+2r^{-i})}\\
=&\sum_{i=1}^n d_i \frac{2r^{-d_i}}{(2^{-ik}+2r^{-i})}\prod_{j=1}^n (2^{-jk}+
2r^{-j})^{d_j}
\end{align*}
Observe that $g(r)\to 0$ as $r\to \infty$ therefore for large $r>0$ we have $N_1 g(r) < \epsilon$. Thus
$\mu(\cup C_{x}^k)-\mu(\cup (\cup E_{y}^r))\leq N_1g(r)<\epsilon$.
This creates a lower bound for the quantity:
$$\mu (A)-2\epsilon
\leq \mu (\cup C^k_{x})-\epsilon\leq
\mu (\cup (\cup E_{y}^r)) \leq \frac{\delta_{1/r}(\mathbb{Z})\cap A}{r^{d'}}
$$

For the upper bound, we cover $U_1\cup \text{int}A$ with $D_{x}^k$,
$\overline{A}\subset U_1\cup \text{int}A$, $\overline{A}$ is closed and bounded, thus compact, so it also has a finite subcover of $D_{x}^k$. We once more assume without loss of
generality that $k>0$ is the same for all of those $D_{x}^k$. We subdivide $\overline{D^k_{x}}$ to
$C^k_{y}$ with foreign interiors. Obviously $A\subset \overline{A}\subset \cup D^k_{x}$ and $\mu(C^k_{y})=\frac{1}{2^{d'k}}$, thus we obtain: 

$$\frac{\#(\delta_{2^{-k}}(\mathbb{Z})\cap A)}{2^{d'k}}\leq \mu(\cup C^{k}_{y})=\mu (\cup D_{x}^k)\leq \mu (A\cup U_1) \leq \mu (A) + \epsilon $$

Once more we want the bound to be valid for all $r>0$. We denote by $N_2:=
\# (\delta_{2^{-k}}(\mathbb{Z})\cap (\cup C_{y}^k))$ the cardinality of the cover. Like previously we will cover it with $E^r_{x}$.
We take the union of all $E_{x}^r$ for which there exist $C^k_{y}$ in the cover such that $E_{x}^r \cap C^k_{y}\neq \emptyset$.
Once more we can observe that 
$\mu(\cup E_{x}^r)-\mu(C_{y}^k)< g(r)$.
Thus 
$\mu (\cup(\cup E_{x}^r))-\mu (\cup C_{y}^k)\leq N_2g(r)$.
If $r>2^k$ large enough then $N_2g(r)<\epsilon$ and we have obtained the upper bound:

$$
\frac{\# (\delta_{1/r}(\mathbb{Z})\cap A)}{r^{d'}}\leq
\mu (\cup (\cup E_{x}^r))\leq
\mu (\cup C_{y}^k)+\epsilon \leq
\mu (U_1\cup A)+\epsilon\leq
\mu (A) +2\epsilon
$$
Therefore $\mu (A) -2\epsilon\leq \frac{\# (\delta_{1/r}(\mathbb{Z})\cap A)}{r^{d'}}\leq \mu (A) +2\epsilon$ and this completes the proof.
\end{proof}

Let $\Gamma < G$ be a lattice in a nilpotent Lie group, $\mathfrak{g}$ its Lie algebra and $\mathfrak{g}_{\infty}=\oplus_{i=1}^{s}V_i$ the graded algebra. Let $S$ be a finite generating set for $\Gamma$ and $d_{CC}$ the Pansu limit metric on $G$ induced by $S$. We will state a series of lemmas that will provide information on how the lattice is rescaled through the dilation inside the stratified group and then we will prove Theorem \ref{Haartheo}. Even though the Haar measure obviously depends on the group structure in this case from both structures arise the same measure up to a scalar multiple. 

We endow $G$ with a strong Malcev basis strongly based on $\Gamma$, recall that $\Psi:\mathbb{R}^d \to G$ are the strong Malcev coordinates for $(G,\ast)$, defined as $\Psi(t_1,\ldots,t_d)=\exp(t_1 X_1)\ast\ldots \ast \exp(t_d X_d)$, and $\Psi_\infty:\mathbb{R}^d \to G$ are the strong Malcev coordinates for $(G,\ast_\infty)$, defined as $\Psi_\infty(t_1,\ldots,t_d)=\exp(t_1 X_1)\ast_\infty\ldots \ast_\infty \exp(t_d X_d)$. The Malcev basis is strongly based on $\Gamma$ with respect to the $(G,\ast)$ structure, i.e. $\Psi(\mathbb{Z}^d)=\Gamma$, but $\Psi_\infty(\mathbb{Z}^d)$ isn't necessary even a group. The exponential coordinates are common for the two group structures, $\Phi:\mathbb{R}^d \to G$ defined as $\Phi(t_1,\ldots,t_d)=\exp(t_1X_1+\ldots t_dX_d)$ and push-forward the Lebesgue measure from $\mathbb{R}^d$ to $G$.

\begin{lemma}\label{Haar}
Let $\nu$ be a Haar measure on $(G,\ast_\infty)$. Then there is a constant $c>0$ such that for every $A\subset G$ bounded $\nu$-continuity set
$$\lim_{r\to\infty}
\frac{\# \left( \delta_{1/r}(\Psi_\infty(\mathbb{Z}^d))\cap A\right)}{r^{d_\Gamma}} = c \mathcal{\nu}
(A).$$
\end{lemma}
\begin{proof}
From Lemma \ref{3.1}, for $B\subset \mathbb{R}^d$ a bounded Lebesgue $\mu$-continuity set, we obtain that $ \# \vert \delta_{1/r}(\mathbb{Z}^d)\cap B \vert / r^{d_\Gamma}\to \mu (B)$. We know that $\Psi_\infty$ pushes-forward the Lebesgue measure to a Haar measure on $(G,\ast_\infty)$. By definition of the dilations we obtain that
\begin{align*}
\# \left(\delta_{1/r}(\Psi_\infty (\mathbb{Z}^d))\cap A \right)&=\#\left(\Psi_\infty (\delta_{1/r}(\mathbb{Z}^d))\cap \Psi_\infty(\Phi^{-1}(A))\right)=\\
&=\#\left(\Psi_\infty\left(\delta_{1/r}(\mathbb{Z}^d)\cap \Psi_\infty^{-1}(A)\right)\right)=\\
&=\#\left(\delta_{1/r}(\mathbb{Z}^d)\cap \Psi_\infty^{-1}(A)\right).
\end{align*}
Therefore, $\frac{1}{r^{d_\Gamma}}\# \vert \delta_{1/r}(\Psi_\infty(\mathbb{Z}^d))\cap A) \vert \to \mu(\Psi_\infty^{-1}(A))$, which is a Haar measure on $(G,\ast_\infty)$. It is known that the Haar measure is unique up to a constant.
\end{proof}

Breuillard and Le Donne \cite{BrL} observe that dilations satisfy $[X,Y]_{\infty}=\lim_{\lambda\to \infty} \delta^{-1}_\lambda [ \delta_\lambda X,\delta_\lambda Y]$ and deduce a new form for the BCH formula. We need some notation before we express the BCH formula. We denote with $s_i$ the number for which $X_i\in V_{s_i}$, let $\vec{s}=(s,s_{d-1},\ldots,1)$ and $\vec{s}\alpha = \sum s_i \alpha_i$, we write $\vec{v}^{\alpha}=v_1^{\alpha_1}\ldots v_d^{\alpha_d}$, for $\alpha\in \mathbb{N}^d$ and $(v)_i=v_i$ for $v=(v_1,\ldots,v_d)\in \mathbb{R}^d$. Recall that $X\bullet Y:=\log(\exp X \ast \exp Y)$ for every $X,Y\in \mathfrak{g}$. For $X=(x_1,\ldots,x_d),Y=(y_1,\ldots,y_d)\in \mathfrak{g}$ the BCH formula is expressed as:
\begin{align}
(X\bullet Y)_i &= x_i + y_i + \sum_{I_i} C_{\alpha,\beta} X^{\alpha}Y^{\beta} + \sum_{J_i} D_{\alpha,\beta} X^{\alpha}Y^{\beta}\label{BHC2}\\
(X\bullet_\infty Y)_i &= x_i + y_i + \sum_{I_i} C_{\alpha,\beta} X^{\alpha}Y^{\beta}\label{BHC3}
\end{align}
where $I_i=\lbrace\alpha , \beta\in\mathbb{Z}^d : \vec{s}\alpha \geq 1 \text{, } \vec{s}\beta\geq 1 \text{ and }
\vec{s}\alpha + \vec{s}\beta = s_i\rbrace$, $J_i=\lbrace\alpha , \beta\in\mathbb{Z}^d : \vec{s}\alpha \geq 1 \text{, } \vec{s}\beta\geq 1 \text{ and }
\vec{s}\alpha + \vec{s}\beta < s_i\rbrace$ 
and $C_{\alpha,\beta}$, $D_{\alpha,\beta}$ are unique constants depending on the Lie structure.

We will enrich this result. We call weighted degree of a monomial $x_1^{\alpha_1}\ldots x_d^{\alpha_d}$ to be $\sum_i \alpha_i s_i$. Additionally, if we write $ x^{(i)}:=(0,\ldots  x_i,\ldots,0)$ for $  x=( x_1,\ldots, x_d)\in \mathbb{R}^d$, then the Malcev coordinates
can be written as $\Psi(  x)=\exp( x^{(1)}\bullet\ldots\bullet  x^{(d)})$ and $\Psi_{\infty}(  x)=\exp ( x^{(1)}\bullet_\infty\ldots\bullet_\infty x^{(d)})$.

\begin{lemma}\label{PsiInf}
There are unique constants $C_\alpha, D_\alpha$ for $\alpha\in\mathbb{N}^d$ such that
\begin{align*}
(\log(\Psi( x)))_i &= x_i   + \sum_{K_i} C_{\alpha } x_{i+1}^{\alpha_{i+1}}\ldots x_{d}^{\alpha_d} + \sum_{L_i} D_{\alpha }  x_{i+1}^{\alpha_{i+1}}\ldots x_{d}^{\alpha_d}\\
(\log (\Psi_\infty( x)))_i &= x_i + \sum_{K_i} C_{\alpha } x_{i+1}^{\alpha_{i+1}}\ldots x_{d}^{\alpha_d}
\end{align*}
for every $ x\in \mathbb{R}^d$, where $K_i = \lbrace \alpha\in \mathbb{N}^d: \vec{s}\alpha= s_i\text{ and } a_j=0 \text{ for } j\leq i\rbrace$ and $L_i = \lbrace \alpha\in \mathbb{N}^d: \vec{s}\alpha < s_i\text{ and } a_j=0 \text{ for } j\leq i\rbrace$.
\end{lemma}
\begin{proof}
We will calculate the product and the constants recursively, applying the BCH formula (\ref{BHC2}) on $ x^{(d-k)}$ and $v_k:= x^{(d-k+1)}\bullet\ldots\bullet  x^{(d)}$. For $ x^{(d-1)}$
and $ x^{(d)}$ the formula is written as 
\begin{align*}
( x^{(d-1)}\bullet  x^{(d)})_i = 
     \begin{cases}
        x_d\text{, } &\quad\text{if }i=d\\
        x_{d-1}\text{, } &\quad\text{if }i=d-1\\ 
       \sum_{I_i} C_{\alpha,\beta}   x_{d-1}^{\alpha_{d-1}} x_d^{\beta_d} + \sum_{J_i} D_{\alpha,\beta}   x_{d-1}^{\alpha_{d-1}} x_d^{\beta_d} &\quad\text{otherwise.} \\ 
     \end{cases}
\end{align*}   
We set $C^{(1)}_{\alpha'}:= C_{\alpha,\beta}$ for $\alpha=(0,\ldots,\alpha_{d-1},0)$, $\beta=(0,\ldots,\beta_d)$ and $\alpha'=(0,\ldots,\alpha_{d-1},\beta_d)$, correspondingly we set $D^{(1)}_{\alpha'}:= D_{\alpha,\beta}$, we set the rest to be equal 0. We calculate recursively, assume that for every $k$ the following relation holds;
\begin{align*}
(v_k)_i = 
     \begin{cases}
         x_i   + \sum_{K_i} C^{(k)}_{\alpha }  x_{i+1}^{\alpha_{i+1}}\ldots x_{d}^{\alpha_d} + \sum_{L_i} D^{(k)}_{\alpha }  x_{i+1}^{\alpha_{i+1}}\ldots x_{d}^{\alpha_d}\text{, } &\text{if }i\geq d-k\\
       \sum_{K_i} C^{(k)}_{\alpha }   x_{d-k}^{\alpha_{d-k}}\ldots x_{d}^{\alpha_d} + \sum_{L_i} D^{(k)}_{\alpha }   x_{d-k}^{\alpha_{d-k}}\ldots x_{d}^{\alpha_d} &\text{otherwise.} \\ 
     \end{cases}
\end{align*} 
and that $C^{(n)}_\alpha$, $D^{(n)}_\alpha$ are defined such that $C^{(n)}_\alpha=D^{(n)}_\alpha=0$, for $j<k$ that $\alpha_j\neq 0$. We will calculate $v_{k+1}$, apply the BCH formula (\ref{BHC2}) on $v_{k+1}= x^{(d-k-1)} \bullet v_k$, where $v_k$ is provided by the previous relation.
\begin{align*} 
(v_{k+1})_i = 
     \begin{cases}
         x_i   + \sum_{K_i} C^{(k)}_{\alpha }  x_{i+1}^{\alpha_{i+1}}\ldots x_{d}^{\alpha_d} + \sum_{L_i} D^{(k)}_{\alpha }  x_{i+1}^{\alpha_{i+1}}\ldots x_{d}^{\alpha_d}\text{, } &\text{if }i\geq d-k-1\\
      (v_k)_i + \sum_{I_i} C_{\alpha,\beta} x_{d-k-1}^{\alpha_{d-k-1}}v_k^{\beta}+\sum_{J_i} D_{\alpha,\beta} x^{\alpha_{d-k-1}}_{d-k-1}v_k^{\beta}
 &\text{otherwise.} \\ 
     \end{cases}
\end{align*} 
For $i\geq d-k-1$ we set $C^{(k+1)}_{\alpha}:=C^{(k)}_{\alpha}$ and $D^{(k+1)}_{\alpha}:=D^{(k)}_{\alpha}$, for the rest we substitute the $v_k^{\beta}$ with the sum and obtain the new constants.

We repeat the procedure for $\Psi_\infty ( x)$. We set $u_k:= x^{(d-k)}\bullet_{\infty}\ldots  x^{(d)}$, and repeatedly apply the BCH formula (\ref{BHC3}). First we apply the BCH formula (\ref{BHC3}) on $x^{(d-1)}$ and $x^{(d)}$.
\begin{align*}   
( x^{(d-1)}\bullet_\infty  x^{(d)})_i = 
     \begin{cases}
        x_d\text{, } &\quad\text{if }i=d\\
        x_{d-1}\text{, } &\quad\text{if }i=d-1\\ 
       \sum_{I_i} C_{\alpha,\beta}   x_{d-1}^{\alpha_{d-1}} x_d^{\beta_d}  &\quad\text{otherwise.} \\ 
     \end{cases}
\end{align*} 
We assume that for $k$ the following relation holds;
\begin{align*}    
(u_k)_i = 
     \begin{cases}
         x_i   + \sum_{K_i} C^{(k)}_{\alpha }  x_{i+1}^{\alpha_{i+1}}\ldots x_{d}^{\alpha_d}\text{, } &\text{if }i\geq d-k\\
       \sum_{K_i} C^{(k)}_{\alpha }   x_{d-k}^{\alpha_{d-k}}\ldots x_{d}^{\alpha_d}  &\text{otherwise.} \\ 
     \end{cases}
\end{align*} 
We will calculate $(u_{k+1})= x^{(d-k-1)}\bullet_\infty u_k$ by applying the BCH formula (\ref{BHC3}) on $x^{(d-k-1)}$ and $u_k$, where $u_k$ is provided by the previous relation.
\begin{align*}      
(u_{k+1})_i = 
     \begin{cases}
         x_i   + \sum_{K_i} C^{(k)}_{\alpha }  x_{i+1}^{\alpha_{i+1}}\ldots x_{d}^{\alpha_d}  \text{, } &\text{if }i\geq d-k-1\\
      (u_k)_i + \sum_{I_i} C_{\alpha,\beta} x_{d-k-1}^{\alpha_{d-n-1}}u_k^{\beta}  
 &\text{otherwise.} \\ 
     \end{cases}
\end{align*} 
For $i\geq d-k-1$ we set $C^{(k+1)}_{\alpha}:=C^{(k)}_{\alpha}$. For the rest we substitute the $u_k^{\beta}$ with the sum and obtain the new constants.

The inductive step induces the same constants $C_a^{(k+1)}$ for both procedures on $\Psi( x)$ and $\Psi_\infty( x)$. The important case is when $i<d-k-1$. We split the terms of $(v_{k+1})_i$ into two parts; the first part will be a sum of monomials of weighted degree exactly $s_i$, and the second will be a sum of polynomials of weighted degree smaller than $s_i$. Obviously for $(v_k)_i= \sum_{K_i} C^{(k)}_{\alpha }   x_{d-k}^{\alpha_{d-k}}\ldots x_{d}^{\alpha_d} + \sum_{L_i} D^{(k)}_{\alpha }   x_{d-k}^{\alpha_{d-k}}\ldots x_{d}^{\alpha_d}$ the separation is clear since for each $\alpha\in K_i$ the weighted degree of $C^{(k)}_{\alpha}  x_{d-k}^{\alpha_{d-k}}\ldots
 x_{d }^{\alpha_{d }}$ equals $\sum s_j \alpha_j=s_i$ for the $C^{(k)}_{\alpha}\neq 0$, while the weighted degree of $D^{(k)}_{\alpha}  x_{d-k}^{\alpha_{d-k}}\ldots
 x_{d }^{\alpha_{d }}$ is smaller than $s_i$.

For the second term we write $v_k^{\beta}=\prod_j (v_k)_j^{\beta_j}=\prod_{j=1}^{d-k-1}(v_k)^{\beta_j}_j \prod_{j=d-k}^{d}(v_k)^{\beta_j}_j$.
We observe that every $(v_k)_j^{\beta_j}$ has a polynomial of weighted degree smaller than $s_j \beta_j$ plus polynomials for which each monomial has weighted degree exactly $s_j \beta_j$. The polynomials that attain weighted degree $s_j \beta_j$ are $ x_j   + \sum_{K_j} C^{(k)}_{\alpha } x_{j+1}^{\alpha_{j+1}}\ldots x_{d}^{\alpha_d}$, for $j\leq d-k$ and $\sum_{K_j} C^{(k)}_{\alpha }   x_{d-k}^{\alpha_{d-k}}\ldots x_{d}^{\alpha_d}$, for $j< d-k$. Thus in the product the monomials that attain the maximum weighted degree $s_i$ arise from the product
\begin{align*}
P_{\beta}( x)=\prod_{j=1}^{d-k-1}( x_j   + \sum_{K_j} C^{(k)}_{\alpha } x_{j+1}^{\alpha_{j+1}}\ldots x_{d}^{\alpha_d})^{\beta_j} \prod_{j=d-k}^{d}(C^{(k)}_{\alpha }   x_{d-k}^{\alpha_{d-k}}\ldots x_{d}^{\alpha_d})^{\beta_j}.
\end{align*}
Therefore, the $C^{(k+1)}_\alpha$ are defined to satisfy
\begin{align*}\sum_{K_i} C_{\alpha}^{(k+1)} x_{d-k-1}^{\alpha_{d-k-1}}\ldots x_d^{\alpha_d} = \sum_{K_i} C^{(k)}_{\alpha }   x_{d-k}^{\alpha_{d-k}}\ldots x_{d}^{\alpha_d}
+\sum_{I_i} C_{\alpha,\beta} x_{d-k-1}^{\alpha_{d-k-1}}
P_{\beta}( x).
\end{align*}

Lastly we observe that the new constants for $(u_{k+1})_i$, for $i<d-k-1$, are defined in each step in order to satisfy the relation
\begin{align*}
\sum_{K_i} C_{\alpha}^{(k+1)} x_{d-k-1}^{\alpha_{d-k-1}}\ldots x_d^{\alpha_d} = \sum_{K_i} C^{(k)}_{\alpha }   x_{d-k}^{\alpha_{d-k}}\ldots x_{d}^{\alpha_d}
+\sum_{I_i} C_{\alpha,\beta} x_{d-k-1}^{\alpha_{d-k-1}}
P_{\beta}( x).
\end{align*}
This is the exact same relation as before, so the constants are the same.
\end{proof}
\begin{corollary}\label{CoroQ}
There exist polynomials $Q_i( x)$ such that $$( x^{(1)}\bullet\ldots \bullet
 x^{(d)})_i=( x^{(1)}\bullet_\infty\ldots  \bullet_\infty x^{(d)})_i + Q_i( x_1,\ldots, x_d)$$ and the weighted degree of $Q_i$ is smaller than $s_i$ for every $i=1,\ldots , d$.
\end{corollary}

The last corollary is useful if we treat $\mathfrak{g}$ as a vector space and work with the $\lVert \cdot \rVert_\infty$ norm. Let $g,h\in G$, then $\Vert \log g - \log h \Vert_\infty=\sup\lbrace \vert (\log g)_i-(\log h)_i\vert : i=1,\ldots d\rbrace$. 
Let $A\subset G$, we will denote with $\mathcal{N}_\epsilon (A) = \lbrace g\in G: \text{dist}(\log (g), \log (A))<\epsilon\rbrace$, where $\text{dist}(x,A)=\inf\lbrace \Vert x-y\Vert_\infty:y\in A\rbrace$ for $x\in \mathfrak{g}$ and $A\subset \mathfrak{g}$.

\begin{proposition}\label{PropAbove}
Let $\nu$ be a Haar measure on $(G,\ast_\infty)$. Then there is a constant $c>0$ such that
\begin{align*}
\lim_{r\to\infty}
\frac{1}{r^{d_\Gamma}}\# \left( \delta_{1/r}(\Gamma)\cap A\right)= c\nu (A),
\end{align*}
for every bounded $\nu$-continuity set $A\subset G$.
\end{proposition}
\begin{proof}

We will use the previous Lemma \ref{PsiInf} and Corollary \ref{CoroQ} to compare $\Psi(\mathbb{Z}^d)$ and $\Psi_\infty(\mathbb{Z}^d)$. Since $\Gamma$ is strongly based on the Malcev base, $\Psi(\mathbb{Z}^d)=\Gamma$. Let $\epsilon>0$. We can limit ourselves in $K$ a closed, bounded, connected subset of $G$ that contains $\mathcal{N}_\epsilon (A)$, since $A$ is bounded.

Let $\Psi_\infty (x)\in \Psi_\infty(\mathbb{Z}^d)$ such that $\delta_{1/r}(\Psi_\infty(x))\in K$, then
\[
\delta_{1/r}(\Psi_\infty (x))=\exp (x_1 r^{-s}X_1)\ast_\infty\ldots\ast_\infty\exp(x_d r^{-1} X_d),
\]
where $x_i\in \mathbb{Z}$. Since $K$ is compact there are $M_i$ such that $\vert x_i r^{-s_i}\vert < M_i$ for every $i=1,\ldots, d$. Denote $t_i:= x_i r^{-s_i}$.By Lemma \ref{PsiInf} for $\Psi_\infty$, we can write
$$
(\log ( \delta_{1/r}(\Psi_\infty (x))) )_i = \sum_{K_i} C_{\alpha} t_i^{\alpha_i}\ldots  t_d^{\alpha_d}\quad$$ 
since $\log (\delta_{1/r}(\Psi_\infty (x)))=\delta_{1/r}(\log (\Psi_\infty(x)))$.
By Lemma \ref{PsiInf}, we can write
$$(\log (\Psi (x) ))_i =\sum_{K_i} C_{\alpha} r^{s_i} t_i^{\alpha_i}\ldots  t_d^{\alpha_d}
+\sum_{J_i} D_{\alpha} r^{\vec{s} \alpha } t_i^{\alpha_i}\ldots  t_d^{\alpha_d}$$
and since 
$\log(\delta_{1/r}(\Psi (x) ))= \delta_{1/r}(\log(\Psi (x) ))$, we obtain
$$(\log(\delta_{1/r}(\Psi (x))))_i = (\log ( \delta_{1/r}(\Psi_\infty (x) ) )_i +\frac{1}{r}
 \sum_{J_i} D_{\alpha} \frac{1}{r^{s_i-\vec{s} \alpha+1 }} t_i^{\alpha_i}\ldots  t_d^{\alpha_d}.
$$

Let $\epsilon_n>0$ be a sequence such that $\epsilon_n$ converges to $0$ as $n\to \infty$.
The sum \linebreak
$\sum_{J_i} D_{\alpha} t_i^{\alpha_i}\ldots  t_d^{\alpha_d}\leq 
\sum_{J_i} \vert D_{\alpha}  \vert M_i^{\alpha_i} \ldots  M_d^{\alpha_d}$. For every $\epsilon_n>0$ there exists a $r_n>0$ large enough such that $\frac{1}{r}\sum_{J_i} \vert D_{\alpha}  \vert M_i^{\alpha_i} \ldots  M_d^{\alpha_d}<\epsilon_n$ for every $r>r_n$. We conclude that $\Vert \log (\delta_{1/r}(\Psi_\infty(x)))-\log(\delta_{1/r}(\Psi(x)))\Vert_\infty<\epsilon_n$. Recall that $\Psi(\mathbb{Z}^d)=\Gamma$ and that $\mathcal{N}_\epsilon ( A ) = \lbrace g\in G : \mathrm{dist}(\log(g),\log(A))<\epsilon \rbrace$. Hence
$$
\#\left(\delta_{1/r}(\Gamma)\cap A\right) \leq
\# \left(\delta_{1/r}(\Psi_\infty(\mathbb{Z}^d))\cap\mathcal{N}_{\epsilon_n}(A)\right)$$
and
\begin{align*}
\# \left(\delta_{1/r}(\Psi_\infty(\mathbb{Z}^d))\cap A \right)&\leq \#\left(\delta_{1/r}(\Gamma)\cap \mathcal{N}_{\epsilon_n}(A) \right)\\
&\leq  \#\left(\delta_{1/r}(\Gamma)\cap A \right) + \#\left(\delta_{1/r}(\Gamma)\cap \mathcal{N}_{\epsilon_n}(\partial A) \right)\\
\end{align*}
By Lemma \ref{Haar} we obtain that $\lim\limits_{r\to\infty} \frac{1}{ r^{d_\Gamma}}\#\left(\delta_{1/r}(\Gamma)\cap A\right) \leq c\nu(\mathcal{N}_{\epsilon_n}(A))$ and\linebreak
$\lim\limits_{ r\to\infty} \frac{1}{r^{d_\Gamma}}\#\left(\delta_{1/r}(\Gamma)\cap A\right)\geq
c\nu(A)-\nu(\mathcal{N}_{\epsilon_n}(\partial(A))$. 

The topology induced by the norm $\Vert \cdot \Vert_\infty$ is the same as the manifolds,
thus $\cap \mathcal{N}_{\epsilon_n}(A) = \bar{A}$ and
$\cap \mathcal{N}_{\epsilon_n}(\partial A)= \partial A$. From assumptions $A$ is $\nu$-continuity therefore $\nu(\partial A)=0$, and thus $\lim\limits_{n\to\infty} \nu(\mathcal{N}_{\epsilon_n} ( A) ) = \nu( \cap \mathcal{N}_{\epsilon_n}(A))=\nu(\bar{A})=\nu(A)$ and \linebreak
$$\lim\limits_{n\to\infty}\nu(\mathcal{N}_{2\epsilon_n}(\partial A))=\mathcal{H}(\cap\mathcal{N}_{2\epsilon_n}(A))=\nu(\partial A)=0.$$ The inequalities above were proven for every $n$, thus taking the limit as $n\to \infty$, we can conclude the desired limit.
\end{proof}

\begin{proof}[Proof of Theorem \ref{Haartheo}]
From the previous Proposition \ref{PropAbove} we know that the limit
$\lim\limits_{r\to\infty}  \# \left( \delta_{1/r}(\Gamma)\cap A\right)/c_S r^{d_\Gamma}=c_S^{-1}c\nu (A)$, where $c_S^{-1}c\nu$ is a Haar measure on $(G,\ast_\infty)$. We write $\mathcal{H}:=c_S^{-1}c\nu$. 

To complete the proof we need to prove that the $\mathcal{H}(B_\infty)=1$. The constant $c_S$ is chosen to be the leading coefficient of the asymptotic growth function in Proposition \ref{4.1.2}. Since $\# \vert \delta_{1/r}(\Gamma) \cap B_{CC}(\id,1) \vert = \# \vert \Gamma\cap B_{CC}(\id,r)\vert$ combined with Proposition \ref{4.1.1}, we obtain that $\#  \vert B_{\lVert \cdot \rVert_S}(\id,n-n^{1-\alpha_S})\vert \leq\# \vert \Gamma\cap B_{CC}(\id,r)\vert \leq \# B_{\lVert \cdot \rVert_S}(\id,n+n^{1-\alpha_S})$. Therefore
\begin{align*}
\lim\limits_{r\to\infty} \frac{ \# \left( \delta_{1/r}(\Gamma)\cap B_\infty\right)}{c_S r^{d_\Gamma}}=1.
\end{align*}
Thus $\mathcal{H}(B_\infty)=1$.
\end{proof}

The push-forward of the Lebesgue measure through the exponential coordinates $\Phi$ is the Haar measure on $(G,\ast)$, the same applies for $(G,\ast_\infty)$, but $(G,\ast_\infty)$ has the same exponential coordinates by definition, therefore a Haar measure of $(G,\ast)$ is Haar measure of $(G,\ast_\infty)$. The Haar measure is unique up to a constant, thus $\mathcal{H}$ is the Haar measure for both $(G,\ast_\infty)$ and $(G,\ast)$. 

\begin{note}
In \cite{BrL} the constant $c_S$ is defined, in the torsion-free case, as the volume of the unit ball $B_{CC}(\id,1)$ with respect to the Haar measure that induces a measure that has volume $1$ on the quotient $G/\Gamma$.
In the previous construction we also obtain a Haar measure, $\mu(\Psi^{-1}(\cdot))$, with volume $1$ on the quotient. This is another way to prove that the constants are equal. Let $\nu$ the induced measure on $G/\Gamma$ and observe that $p\circ \Phi\vert_{[0,1)^d}:[0,1)^d\to G\to G/\Gamma$ is a bijection. Thus $\nu^{\ast}(G/\Gamma)=\mu([0,1)^d)=1$.
\end{note}

\section{Background on inhomogeneous random graphs}\label{4}
\subsection{Percolation}
We now mention a celebrated result on percolation on graphs by Liggett--Schonmann--Stacey \cite{Lig}, which will be crucial for the proof of the main theorem. Given $n\in\N$, a percolation measure on a graph $G$ is called \emph{$n$-independent} if for any two sets $E_1,E_2$ of edges at graph distance at least $n$ the value on each edge of $E_1$ is independent of each value of the edges in $E_2$. The following proposition is a special case of their result. We consider the infinite integer lattice $\mathbb{Z}^2$.

\begin{proposition}\label{Lig}
For every $K>0$ there exists a $p_K<1$ with the property that for any $K$-independent bond percolation on $\mathbb{Z}^2$ that the probability for an edge to be open is bigger than $p_K$, there is an infinite connected component of open edges in $\mathbb{Z}^2$ almost surely.
\end{proposition}

\subsection{Inhomogeneous random graphs} 
The rest of this section is introduction to terminology and a short description of results on inhomogeneous random graphs derived from \cite{BJR2}.
The notions are adapted to our context, thus they will be stated less generally to avoid unnecessary complexity. The reason we state the following results is because we are going to study the spread-out percolation locally and treat it as an inhomogeneous random graph.

Let $\mathcal{S}$ be a separable metric space 
endowed with $\mu$ a Borel measure such that $0<\mu(\mathcal{S})< \infty$.
For every $\rho\in (0,\infty)$, let $V_\rho$ be a collection of points in $\mathcal{S}$. The triple
$(\mathcal{S},\mu,(V_\rho )_{\rho>0} )$ is called a \emph{generalized vertex space} if $\lim_{\rho}\#(V_\rho \cap A)/\rho =\mu (A)$ for every $\mu$-continuity set $A$, where a set $A\subseteq \mathcal{S}$ is called \emph{$\mu$-continuity} if it is measurable and $\mu(\partial A)=0$. The boundary $\partial A$ is defined as the elements in the closure of $A$ that do not belong in the interior of $A$.

We call a symmetric, non-negative, measurable function $\kappa:\mathcal{S}\times \mathcal{S}\to \mathbb{R}$  a \emph{kernel}. A kernel $\kappa$ is called \emph{irreducible} if, for every set $A\subseteq \mathcal{S}$, we have that that $\kappa\equiv 0$ almost everywhere on $A\times (\mathcal{S}\setminus A)$ implies either $\mu (A)=0$ or $\mu(\mathcal{S}\setminus A)=0$.

\begin{definition}\label{2.inho}
We form the inhomogeneous random graph $G(\rho,\kappa)$, this is the random graph with vertex set $V_\rho$, for each two vertices $x,y\in V_\rho$ we connect them independently with probability $\min \lbrace \frac{\kappa(x,y)}{\rho},1\rbrace$.
\end{definition}

For every graph $G$ we denote by $e(G)$ the number of edges in the graph and by $C_i(G)$ the number of edges in the $i$th largest connected component of the graph. If there is no such a component then $C_i(G)=0$. 

A sequence of kernels $(\kappa_\rho)$ on $\mathcal{S}$ is called graphical on $(\mathcal{S},\mu,(V_\rho)_\rho))$ with limit $\kappa$ if $x_\rho\to x$ and $y_\rho\to y$ implies $\kappa_\rho(x_\rho,y_\rho)\to \kappa (x,y)$ for almost every $(x,y)\in \mathcal{S}\times\mathcal{S}$ and $\frac{1}{\rho}\mathbb{E}\left[ e(G(\rho,\kappa_\rho))\right]\to\frac{1}{2}
\iint\limits_{\mathcal{S}^2} \kappa(x,y) \mu(\mathrm{d} x) \mu(\mathrm{d}y)$, where $\kappa$ is a kernel on $(\mathcal{S},\mu,(V_\rho)_\rho)$, continuous almost everywhere and $\kappa \in L^1$.

For every kernel $\kappa$ on $(\mathcal{S},\mu)$ an integral operator $T_\kappa$ is defined, for every $f:\mathcal{S}\to \mathbb{R}$ measurable function in $L^2(\mathcal{S}\times\mathcal{S},\mu\times\mu)$ set $T_\kappa f (x) = \int_{\mathcal{S}} \kappa (x,y) f(y) \mu (\mathrm{d} y)$. Correspondingly its norm is defined as $\Vert T_\kappa\Vert =\sup\lbrace \Vert T_\kappa f\Vert_2 :
f\geq 0 \text{ measurable}, \Vert f\Vert_2 \leq 1\rbrace$.
\begin{proposition}[{\cite[Theorem 3.1]{BJR2}}]\label{BolProp1}
Let $(\kappa_{\rho})$ be a graphical sequence of irreducible kernels on $(\mathcal{S},\mu,(V_\rho))$ with limit $\kappa$, such that $\Vert T_\kappa\Vert >1$, then there exists an $\alpha(\kappa)>0$ such that $C_1(G(\rho,\kappa_{\rho}))/\rho$ converges in probability to $\alpha(\kappa)$.
\end{proposition}
\begin{proposition}[{\cite[Theorem 3.6]{BJR2}}]\label{BolProp2}
Let $(\kappa_{\rho})$ be a graphical sequence of irreducible kernels on $(\mathcal{S},\mu,(V_\rho))$ with irreducible limit $\kappa$, then $C_2(G(\rho,\kappa_{\rho}))=o_p({\rho})$, i.e. that $C_2(G(\rho,\kappa_{\rho}))/\rho$ converges in probability to $0$.
\end{proposition}

\section{The case of torsion-free nilpotent groups}\label{5}
We now prove Theorem \ref{4.1}. 
In order to study the spread-out percolation in the group, we proceed by defining a similar random graph but qualitatively different and compare those two. 

Let $\Gamma < G$ be a lattice in a nilpotent Lie group that it is not isomorphic to $\mathbb{Z}$. Let $\mathfrak{g}$ be the Lie algebra of $G$ and $\mathfrak{g}_{\infty}=\oplus_{i=1}^{s}V_i$ the graded algebra. Let $S$ be a finite generating set for $\Gamma$ and $d_{CC}$ the Pansu limit metric on $G$ induced by $S$.

For every degree parameter $\lambda>0$ and distance parameter $r>0$ we form $G^{r,\lambda}:=G^{r,\lambda}(G,\Gamma,S)$, the random graph that has for vertex set $\Gamma$, for each pair of vertices $\gamma_1,\gamma_2\in\Gamma$ at distance $d_{CC} (\id, \gamma_2^{-1}\ast \gamma_1) \leq r$ we connect them independently by an edge with probability $\frac{\lambda}{c_S r^{d_\Gamma}}$. 

Recall that we denote as $B_\infty:=B_{CC}(\id,1)$, then $B_{CC}(\id,r)=\delta_r(B_\infty)$ for $r>0$. Note that $G^{r,\lambda}$, likewise $G_{r,\lambda}$, depends on the set of generators $S$, because $B_\infty$ is the unit ball with respect to the Pansu limit metric $d_{CC}$, which is induced by $S$. This graph has the advantage though that the distribution, based on which it is defined, has nice asymptotic properties. We can express the probability that
we connect any pair of vertices $\gamma_1,\gamma_2\in \Gamma$ as $\min\lbrace 1,
\frac{\lambda}{c_S r^{d_\Gamma}}\mathds{1}_{\delta_r(B_\infty)}(\gamma_2^{-1}\ast \gamma_1)\rbrace$.

\begin{restate}
For every $\lambda >1$ there exists $R>0$ such that if $r>R$, then $G^{r,\lambda}$ has an infinite component with probability 1.
\end{restate}

Before we move on to studying $G^{r,\lambda}$ and proving this theorem, we will show how it implies Theorem \ref{4.1}.

\begin{proof}[Proof of Theorem \ref{4.1}]
For a degree parameter $\lambda'\in (1,\lambda)$ and a suitable choice of distance parameter $r'>0$, we will show that the probability two vertices in $G^{r,\lambda'}$ are connected by an edge is dominated by the probability that the same vertices are connected with an edge in $G_{r',\lambda}(\Gamma)$. Combined with the previous theorem, the validation of this statement completes the proof.

Crucial for the comparison of the two random graphs are the results on comparison of metrics on $\Gamma$ and the bound for the growth function of $B_{\lVert \cdot \rVert_S}$. From (\ref{4.1.1}) there is a function $f(r)=O(r^\alpha)$, where $\alpha<1$, such that $\delta_r(B_\infty)\cap \Gamma \subset B_{\lVert \cdot \rVert_{S}}(\id,r+f(r))$ for every $r>0$. From (\ref{4.1.2}), there is a function $g(r)=O(r^\beta)$, where $\beta<d_\Gamma$ such that $\# (B_{\lVert \cdot \rVert_S}(\id,r))\geq c_S r^{d_\Gamma}+g(r)$.
Combining those two we obtain
\begin{align*}
\# (B_{\lVert \cdot \rVert_S}(\id,r+f(r)))&\leq c_S (r+f(r))^{d_\Gamma}+g(r+f(r))\\
&\leq c_S r^{d_\Gamma} + h(r)
\end{align*}
where $h(r)=o(r^{d_\Gamma})$.

There exists an $r_1\in \mathbb{N}$ large enough such that for every $r>r_1$
$$\frac{\lambda'}{c_S r^{d_\Gamma}}\mathds{1}_{\delta_r (B_{\infty})} \leq 
\frac{\lambda}{\# B_{\lVert \cdot \rVert_S}(\id,r+f(r))}\mathds{1}_{ B_{\lVert \cdot \rVert_S}(\id,r+f(r))}$$
Therefore the probability that two vertices in $G^{r,\lambda'}(\Gamma)$ are connected by an edge is dominated by the probability those exact two vertices are connected by an edge in $G_{r+f(r),\lambda}(\Gamma)$. From Theorem \ref{4.2}, there exists an $r_2$ such that for any $r>r_2$ the graph $G^{r,\lambda'}(\Gamma)$ has an infinite component almost surely. We conclude that for $r>\max\lbrace r_1+f(r_1),r_2+f(r_2)\rbrace$ the $G_{r,\lambda}(\Gamma)$ has an infinite connected component almost surely.
\end{proof}

For the proof of Theorem \ref{4.2} we will study the graph $G^{r,\lambda}$ locally by shrinking the graph through dilations and focusing on the finite part that belongs in a fixed region $\mathcal{S}\subset G$, to be specified later. 
The random graph $\delta_{1/r} ( G^{r,\lambda})\cap \mathcal{S}$ has the same distribution as the graph that has for vertices $\delta_{1/r}(\Gamma)\cap \mathcal{S}$ and each pair of vertices $g,h\in   \delta_{1/r}(\Gamma)\cap \mathcal{S}$ is connected independently by an edge with probability $\min\lbrace 1, \frac{\lambda}{c_S r^{d_\Gamma}}\mathds{1}_{B_{\infty}}(\delta_{1/r}(\delta_r(h^{-1})\ast \delta_r(g)))\rbrace$.
Expressing the graph in this form is more complicated but it allows us to treat it as an inhomogeneous random graph $G(\rho,\kappa_\rho)$ as in Definition \ref{2.inho} in the set up of the next proposition. 

\begin{proposition}\label{Prop1}
Let $\mathcal{S}$ be a bounded, open, connected subset of $G$. We set $\rho:= c_S r^{d_\Gamma}$ and we consider the sequence of vertex sets $V_\rho := \delta_{1/r}(\Gamma)\cap \mathcal{S}$. We define $\kappa: \mathcal{S} \times  \mathcal{S} \to \mathbb{R}$ as $\kappa:=\lambda \mathds{1}_{B\infty}(h^{-1}\ast_\infty g)$ and $\kappa_\rho: \mathcal{S}\times  \mathcal{S} \to \mathbb{R}$ as $\kappa_\rho(g,h):=\lambda\mathds{1}_{B_{\infty}}(\delta_{1/r}(\delta_r(h^{-1})\ast \delta_r(g)))$.
Then the triple $(\mathcal{S},\mathcal{H},(V_\rho)_{\rho})$ is a generalized vertex space, $\kappa$ is an irreducible kernel and $(\kappa_\rho)_{\rho>0}$ is a sequence of kernels graphical on $(\mathcal{S},\mathcal{H},(V_\rho)_{\rho})$ with limit $\kappa$.
\end{proposition}

We will need some lemmas for the proof of Proposition \ref{Prop1}, their proofs are technical, thus they will be given in the end of this section. We define $F_r:=\delta_{1/r} \circ
\Psi \circ \delta_r$.
\begin{lemma}\label{Psilim}
For every $x\in\mathbb{Z}^d$ the sequence $F_r(x)\to \Psi_\infty (x)$ as $r\to \infty$. Moreover the convergence is uniform on bounded sets.
\end{lemma}

\begin{lemma}\label{Bounded}
Let $\mathcal{S}$ be a bounded, closed, connected subset of $G$. For every $\eta >0$ there is an $r_0>0$ such that for $r>r_0$ 
$$
\Psi_\infty^{-1}(\mathcal{S}\setminus \mathcal{N}_\eta (\partial \mathcal{S}))\subseteq F_r^{-1}(\mathcal{S})\subseteq \Psi^{-1}_\infty(\mathcal{N}_\eta(\mathcal{S})).
$$
Moreover there exists a $B\subset \mathbb{R}^d$ such that $F_r^{-1}(\mathcal{S})\subset B$ for every $r>0$. 
\end{lemma}

\begin{lemma}\label{ExtraLe}
Let $\rho=c_Sr^{d_\Gamma}$ and $(g_n),(h_n)$ two sequences in $G$ such that $g_n\to g$ and $h_n\to h$, as $n \to \infty$. Then $\delta_{1/r}(\delta_r (g_\rho)\ast \delta_r(h_\rho))\to g\ast_\infty h$, as $\rho \to \infty$.
\end{lemma}

\begin{proof}[Proof of Proposition \ref{Prop1}]
We begin by observing that the triple $(\mathcal{S},\mathcal{H},(V_\rho)_{\rho})$ is indeed a generalized vertex space, it is an immediate corollary of Theorem \ref{Haartheo} since $\mathcal{H}(\mathcal{S})<\infty$. Now we need to prove a series of properties for $\kappa$ and $\kappa_\rho$.

Firstly, $\kappa_{\rho}$ and $\kappa$ are kernels. Indeed, they are positive functions. The fact that they are symmetric arises from the inverse element in the group structure, which coincides in the two different group structures. This can be derived from the BCH formula (\ref{BHC}) for the Lie brackets. Since the groups are Lie groups the inverse element as a function is continuous. Both of them are positive functions, thus they are indeed kernels. 

Moreover, $\kappa$ is continuous a.e. and  $\kappa\in L^1(\mathcal{S}\times\mathcal{S},\mathcal{H}\otimes\mathcal{H})$. This is a consequence of $\mathcal{H}(\partial B_{\infty})=0$, which is valid because $\partial B_\infty \subset B_{CC}(\id,1+\epsilon)\setminus B_{CC}(\id,1-\epsilon)$ and $\mathcal{H}(B_{CC}(\id,r))=r^{d_\Gamma}$. Thus $\kappa$ is almost everywhere continuous.
Additionally, $\Vert \kappa\Vert_1 = \iint\limits_{\mathcal{S}\times\mathcal{S}} \kappa(x,y) \mathcal{H}\otimes\mathcal{H}(\mathrm{d} (x,y))\leq \lambda \iint\limits_{\mathcal{S}\times\mathcal{S}} \mathcal{H}\otimes\mathcal{H}(\mathrm{d} (x,y))\leq \lambda \mathcal{H}(\mathcal{S})^2\leq\infty$ as required.

Next, $\kappa$ is an irreducible kernel.
Let $A\subset \mathcal{S}$ such that $0<\mathcal{H}(A)<\mathcal{H}(\mathcal{S})$. Then there exists $g\in G$ such that $B_{CC}(g,\frac{1}{2})\cap A$ and $B_{CC}(g,\frac{1}{2})\cap A^c$ both have positive measures, because $\mathcal{S}$ is connected and open. 
Indeed, there are two balls $B_1,B_2$ such that $B_1\cap A$ and $B_2\cap A^{c}$ have positive measures, we connect their centers with a line and cover it with open balls, since it is a compact set, it has a finite subcover. Thus we have a finite number of balls intersecting each other, if there is no ball intersecting both $A$ and $A^{c}$ in a set of positive measure then there are two balls with non empty intersection such that one intersects $A$ in a set of positive measure and the other intersects $A^{c}$ in a set of positive measure, yet their intersection has positive measure, which is a contradiction. 
We obtain $\kappa(g_1,g_2)\neq 0$ for every $(g_1,g_2)\in B_{CC}(g,\frac{1}{2})\cap A\times B_{CC}(g,\frac{1}{2})\cap A^c$, since $g^{-1}_2\ast_\infty g_1\in B_\infty$ for any $g_1,g_2\in B_{CC}(g,\frac{1}{2})$, which concludes that $\kappa$ is irreducible.

We move on to prove that $\kappa_\rho$ is a graphical sequence with limit $\kappa$. Thus, we need to prove that if $g_\rho \to g$ and $h_\rho \to h$, as $\rho \to \infty$ then $\kappa_\rho(g_\rho,h_\rho)\to \kappa(g,h)$, as $\rho\to \infty$, almost everywhere in $G\times G$. 
Let $g_\rho,h_\rho\in G$ such that $g_\rho\to g$ and $h_\rho\to h$. It is enough to prove that $\kappa_\rho(g_\rho,h_\rho^{-1})\to \kappa(g,h^{-1})$ because the inverse element as a function is continuous in the Lie group. 
Lemma \ref{ExtraLe} tells us that $\delta_{1/r}(\delta_r (g_\rho)\ast \delta_r(h_\rho))\to g\ast_\infty h$.

We recall that $\kappa$ is defined as $\lambda\mathds{1}_{B_\infty}(h^{-1}\ast_\infty g)$.
Thus the kernel is discontinuous on points of the boundary, $h^{-1}\ast_\infty g \in \partial B_\infty$. We define $U:=\lbrace (x,y)\in G\times G : y^{-1}\ast_\infty x\in \partial B_\infty \rbrace$, $U$ is a closed set.
We treat $\mathfrak{g}$ as $\mathbb{R}^d$ and consider the $\Vert \cdot \Vert_\infty$ norm on it, correspondingly for $X\in\mathfrak{g}$ and $A\subset \mathfrak{g}$ the $\mathrm{dist}(X,A)=\inf\lbrace\Vert X-Y\Vert :Y\in A\rbrace$.
For any $(g,h)\in U^c$ there exists an $\epsilon>0$ such that $\text{dist} (\log(h^{-1}\ast_\infty g), \log(\partial B_\infty)) > \epsilon$, because the boundary $\partial B_\infty$ is a closed set. Since the sequence $\delta_{1/r}(\delta_r(h^{-1}_\rho)\ast\delta_r(g_\rho))$ converges to $h^{-1}\ast_\infty g$, therefore there exists an $R>0$ such that for every $r>R$ the $\Vert \log (\delta_{1/r}(\delta_r(h^{-1}_\rho)\ast\delta_r(g_\rho)))-\log(h^{-1}\ast_\infty g)\Vert_\infty \leq \epsilon$. Thus $\kappa_\rho(g_\rho,h_\rho)=\kappa(g,h)$ for every $r>R$. 

To complete the proof of this convergence we need to prove that the kernel is discontinuous on a set of measure $0$, i.e. $\mathcal{H}\otimes\mathcal{H}(U)=0$, but $U_x=\lbrace y\in G : y^{-1}\ast_\infty x \in \partial B_\infty\rbrace=x\ast_\infty \partial B_\infty$ and $\mathcal{H}(x\ast_\infty \partial B_\infty)=0$, because the boundary of $B_\infty$ has measure $0$. We conclude that almost everywhere if $g_\rho \to g$ and $h_\rho \to h$, then $\kappa_\rho (g_\rho,h_\rho)\to \kappa(g,h)$. 
This concludes the proof that $\kappa_\rho$ is a graphical sequence with limit $\kappa$.

The last property we need to prove is that $\frac{1}{\rho}\mathbb{E}\left[ e(G(\rho,\kappa_\rho))\right] \to \iint\limits_{\mathcal{S}\times\mathcal{S}} \kappa(x,y)\mathcal{H}\otimes\mathcal{H}(\mathrm{d} (x,y))$ as $\rho \to \infty$, where $G(\rho,\kappa_\rho)$ is the inhomogeneous random graph arising by the kernel $\kappa_\rho$ on the generalized vertex space $(\mathcal{S},\mathcal{H}, (V_\rho)_\rho)$.
We want to calculate the limit, therefore we can restrict to $\rho > \lambda$ and substitute 
\begin{align*}
\frac{1}{\rho}\mathbb{E}\left[ e(G(\kappa_\rho,\rho))\right]&=\frac{1}{2\rho}
\sum\limits_{x,y\in V_\rho}\min\left\lbrace \frac{\kappa_\rho(x,y)}{\rho},1\right\rbrace\\
&=\frac{1}{2\rho}
\sum\limits_{x,y\in V_\rho} \frac{\kappa_\rho(x,y)}{\rho}
\end{align*}
The sum is over elements of $V_\rho$ which is a dynamical set, it varies as $\rho$ grows and this causes problems in calculations, what we expound below is done in order to slide over this problem with the help of the previous lemmas. We substitute in the sum $v=F_r^{-1}(x)$ and $u=F_r^{-1}(y)$ and the sum is written as
$$
\sum\limits_{x,y\in V_\rho}\frac{1}{\rho^2} \kappa_\rho(x,y)  = 
\sum\limits_{v,u\in F_r^{-1}(V_\rho)} \frac{1}{\rho^2} \kappa_\rho(F_r(v),F_r(u)).
$$
Recall that $V_\rho=\delta_{1/r}(\Gamma)\cap \mathcal{S}$, $\Psi^{-1}(\Gamma)=\mathbb{Z}^d$ and $\delta_{1/r}(\delta_r(x))=x$, simply substitutions lead to the equality:
\begin{align*}
F^{-1}_r (V_\rho)&=F_r^{-1}(\delta_{1/r}(\Gamma))\cap F^{-1}_r(\mathcal{S})\\
&=\delta_{1/r}(\Psi^{-1}(\delta_r(\delta_{1/r}(\Gamma)) ))\cap F^{-1}_r(\mathcal{S})\\
&=\delta_{1/r}(\mathbb{Z}^d)\cap F_r^{-1}(\mathcal{S}).
\end{align*} 

Lemma \ref{Bounded} allows us to sandwich $F^{-1}_r(\mathcal{S})$, and hence sandwich $F^{-1}_r(V_\rho)$ by $C_r^{\eta}:=\Psi_\infty^{-1} (\mathcal{S}\setminus \mathcal{N}_{\eta}(\partial\mathcal{S}))\cap\delta_{1/r}(\mathbb{Z}^d)$ and $D_r^{\eta}:=\Psi_\infty^{-1} (\mathcal{N}_{\eta}(\mathcal{S}))\cap\delta_{1/r}(\mathbb{Z}^d)$
so that 
$ C_r^{\eta}\subset F_r^{-1}(V_\rho)\subset D_r^{\eta}$. Lastly we set as $A_r:= \Psi_\infty^{-1}(\mathcal{S})\cap \delta_{1/r}(\mathbb{Z}^d)$ and bound the following difference;

\begin{align*}
\left\vert 
 \sum\limits_{v,u\in A_r}\frac{\kappa_\rho(F_r(v),F_r(u))}{(c_S r^{d_\Gamma})^2}
 -\sum\limits_{v,u\in F_r^{-1}(V_\rho)} \frac{\kappa_\rho(F_r(v),F_r(u))}{{(c_S r^{d_\Gamma})^2}}\right\vert  & \leq \sum\limits_{v,u\in D_r^{\eta} \setminus C_r^{\eta}} \frac{ \kappa_\rho(F_r(v),F_r(u))}{{(c_S r^{d_\Gamma})^2}} \\
&\leq  \sum\limits_{v,u\in D_r^{\eta} \setminus C_r^{\eta}}\frac{\lambda}{(c_S r^{d_\Gamma})^2}.
\end{align*}
From Theorem \ref{Haartheo} though we obtain that the right part of the inequality
\[
\lim\limits_{r\to\infty}  \sum\limits_{v,u\in D_r^{\eta} \setminus C_r^{\eta}}\frac{\lambda}{(c_S r^{d_\Gamma})^2} =\lambda \mathcal{H}(\mathcal{N}_{\eta}(\partial \mathcal{S}))^2
\]
for every $\eta>0$, but $\mathcal{H}(\mathcal{N}_{\eta}(\partial \mathcal{S}))\to \mathcal{H}(\partial\mathcal{S})=0$ as $\eta \to 0$.
We deduce the proof to the calculation of the following limit;
\begin{align}\label{3giaprotasi}
\lim_{r\to\infty} \sum\limits_{v,u\in A_r}\kappa_\rho(F_r(v),F_r(u))\frac{1}{(c_S r^{d_\Gamma})^2}=\lim_{r\to\infty} \sum\limits_{v,u\in F_r^{-1}(V_\rho)} \frac{\kappa_\rho(F_r(v),F_r(u))}{{(c_S r^{d_\Gamma})^2}}.
\end{align}

From Lemma \ref{Psilim}, the sequence $F_r(v)$ converges to $\Psi_\infty(v)$ as $r \to \infty$ uniformly on $\Psi_\infty^{-1}(\mathcal{S})$. In combination with the previous property of the sequence of kernels, that if $g_\rho\to g$ and $h_\rho\to h$, then $\kappa_\rho(g_\rho,h_\rho)\to\kappa(g,h)$. We obtain that $\kappa_\rho(F_r(v),F_r(u)) \to \kappa (\Psi_\infty(v),\Psi_\infty(u)) $ as $r \to \infty$ almost everywhere.

Now that we established $\mathcal{A}_r$ a fixed domain to work on and $\kappa_\rho(F_r(v),F_r(u))$ converges pointwise, we will obtain the desired limit from Lebesgue's theorem. Essentially we create a partition of $\Psi_\infty^{-1}(\mathcal{S})\times \Psi_\infty^{-1}(\mathcal{S})$. We define $P^{(r)}_{(w,z)}=E^r_w\times E^r_z$ for every $z,w\in \delta_{1/r}(\mathbb{Z}^d)$, where $E^r_z=[z_1,z_1+r^{-s})\times\ldots[z_d,z_d+r^{-1})$
as in Lemma \ref{3.1}. We define the step functions $f_r(v,u):=\kappa_\rho(F_r(w),F_r(z))$ for every $(v,u)\in P^{(r)}_{(w,z)}$.

We will show that $f_r(v,u)\to \kappa (\Psi_\infty(v),\Psi_\infty(u))$ pointwise, as $r\to \infty$, almost everywhere.
Recall that convergence is not guaranteed for $(x,y)\in U$, where we have defined $U$ as $\lbrace (x,y)\in G\times G:y^{-1}\ast_\infty x \in \partial B_\infty\rbrace$ a set of measure $0$.
Let $(v,u)\in \Psi_\infty^{-1}(\mathcal{S})\times \Psi_\infty^{-1}(\mathcal{S})$ such that $(\Psi_\infty(v),\Psi_\infty(u))\in U^c$, this means that $\Psi_\infty(u)^{-1}\ast_\infty\Psi(v)\in \text{int}(B_\infty)\cup \text{int}(B^c_\infty)$.
Let $g:=\Psi_\infty(u)^{-1}\ast_\infty\Psi(v)$, there exists an open ball $\mathcal{N}_\epsilon(g)$ with respect to $\Vert \cdot \Vert_\infty$ in $G$ such that $\mathcal{N}_\epsilon(g)\subset \text{int}(B_\infty)\cup \text{int}(B^c_\infty)$. For this $\epsilon>0$ we will find $r_0>0$ such that for every $r>r_0$ the following holds; if $(v,u)\in P^{(r)}_{(w,z)}$, then every $(x,y)\in P^{(r)}_{(w,z)}$ has the property $\Psi_\infty(y)^{-1}\ast_\infty\Psi_\infty(x)\in \mathcal{N}_\epsilon(g)$. 
Indeed, $\Psi^{-1}_\infty(\mathcal{S})$ is a bounded set. The coordinates of $\log ( \Psi_\infty(y)^{-1}\ast_\infty\Psi_\infty(x))$ are given by polynomials. Since $(x,y)\in P_{(z,w)}$ and $(v,u)\in P_{(z,w)}$, then $(x)_i-(v)_i < r^{-s_i}$ and $(y)_i-(u)_i < r^{-s_i}$. Therefore if $r_0=\epsilon^{-1}$, then for $r>r_0$ and $(x,y)\in P_{(z,w)}$, 
$\Vert \log ( \Psi_\infty(y)^{-1}\ast_\infty\Psi_\infty(x))- \log ( \Psi_\infty(u)^{-1}\ast_\infty\Psi_\infty(v))\Vert_\infty< \epsilon$.

We conclude that we have shown $f_r(v,u)\to \kappa (\Psi_\infty(v),\Psi_\infty(u))$ pointwise almost everywhere, thus
\begin{align}\label{1giaprotasi}
\iint\limits_{\Psi^{-1}_\infty (\mathcal{S})\times\Psi^{-1}_\infty (\mathcal{S})} f_r(v,u) \mu\otimes\mu (\mathrm{d} (v,u))
\to 
\iint\limits_{\Psi^{-1}_\infty (\mathcal{S})\times\Psi^{-1}_\infty (\mathcal{S})} \kappa (\Psi_\infty(v),\Psi_\infty(u)) \mu\otimes\mu(\mathrm{d} (v,u)).
\end{align}
The first integral evaluates as
\begin{align}
\iint\limits_{\Psi^{-1}_\infty (\mathcal{S})\times\Psi^{-1}_\infty (\mathcal{S})} f_r(v,u) \mu\otimes\mu(\mathrm{d} (v,u))&=\sum_{w,z\in A_r}
\iint\limits_{P_{(w,z)}} f_r(v,u) \mu\otimes\mu(\mathrm{d} (v,u)) \notag \\
&= \sum_{w,z\in A_r}
\iint\limits_{P_{(w,z)}} \kappa_\rho(F_r(w),F_r(z)) \mu\otimes\mu(\mathrm{d} (v,u)) \notag \\
&=\sum_{w,z\in A_r} \kappa_\rho(F_r(w),F_r(z)) \mu\otimes\mu (P_{(w,z)}) \notag \\
&=\sum_{w,z\in A_r} \kappa_\rho(F_r(w),F_r(z)) \frac{1}{(r^{d_\Gamma})^2}.\label{2giaprotasi}
\end{align}
Recall that $\mathcal{H}(\cdot)=\frac{\mu ( \Psi^{-1}_\infty(\cdot))}{c_S}$ and combine (\ref{3giaprotasi}),(\ref{1giaprotasi}) and (\ref{2giaprotasi}) to obtain the desired limit. 
\end{proof}

For every choice of $\mathcal{S}$ an inhomogeneous random graph $G(\rho,\kappa_\rho)$ is formed as in Definition \ref{2.inho}. This random graph will have the same distribution as $\delta_{1/r}(G^{r,\lambda})\cap \mathcal{S}$. Therefore our approach to studying $G^{r,\lambda}$ will be as follows; we will find a specific region $\mathcal{S}$. We will use the tools from the previous construction to describe the graph inside the region $\mathcal{S}$. We will define an action of $\mathbb{Z}^2$ that acts freely on $\mathcal{S}$ and covers a part of $G$. Inside that part we will find the infinite component.

\begin{proof}[Proof of Theorem \ref{4.2}]
We consider the random graph $G'(\Gamma)$, that has vertex set $\delta_{1/r}(\Gamma)$ and for two vertices $\gamma_1,\gamma_2\in\delta_{1/r}(\Gamma)$ we connect them with probability \linebreak $\min\lbrace 1, \frac{1}{cr^{d_\Gamma}}\mathds{1}_{B_\infty}(\delta_{1/r}(\delta_r(\gamma_2)^{-1}\ast\delta_r(\gamma_1))\rbrace$. It is obvious that this graph $G'(\Gamma)$ has the exact same distribution as $G^{r,\lambda}$ rescaled through $\delta_{1/r}$. This is due to the fact that the indicator function can be written as $\mathds{1}_{\delta_r(B_\infty)}(x)=\mathds{1}_{B_\infty}(\delta_{1/r}(x))$ and that for both group structures the inverse element function is the same, $\delta_r (x^{-1})=\delta_r(x)^{-1}$.

We begin with the observation that if $g=\exp (\lambda_1 X_1 )\in Z(G)$, then $g\ast h= g\ast_\infty h$,
for any $h\in G$. Indeed if $Y=(y_1,\ldots,y_d)\in\mathfrak{g}$ such that $\exp Y=h$ then the BCH formula (\ref{BHC2}) and (\ref{BHC3}) calculates the coordinates of $\log (g\ast h)$ and $\log (g\ast_\infty h)$ correspondingly,
$$(\log (g\ast h))_1=y_1+\lambda  \text{ and } (\log (g\ast_\infty h))_1=y_1+\lambda$$
and
$$(\log (g\ast h))_j=y_j+0 \text{ and } (\log (g\ast_\infty h))_j=y_j+0,\quad \text{for every } j>1,
$$
so that $g\ast h=g\ast_\infty h$. The same applies for $g=\exp (\lambda X_2)$, i.e. $g\ast h = g\ast_\infty h$ for any $h\in G$. Indeed, we can assume without loss of generality that $X_2\notin \mathfrak{g}^{(s)}$, otherwise it is the same case as the previous one, $X_1$. 
Let $Y=(y_1,\ldots,y_d)\in \mathfrak{g}$ such that $\exp (Y)=h$, then the BCH formula (\ref{BHC2}) and (\ref{BHC3}) applied for $\ast$ and $\ast_\infty$ correspondingly on the first, second and $j$-th coordinates, provides us with
\begin{align*}
(\log (g\ast h))_1 &=y_1+\sum_{I_1} C_{\alpha,\beta} (\lambda X_2)^\alpha Y^{\beta} + \sum_{J_1} D_{\alpha,\beta} (\lambda X_2)^\alpha Y^{\beta},\\
(\log (g\ast_\infty h))_1 &=y_1+\sum_{I_1} C_{\alpha,\beta} (\lambda X_2)^\alpha Y^{\beta},\\
(\log (g\ast h))_2 &=y_2+\lambda  \text{ and } (\log (g\ast_\infty h))_2=y_2+\lambda,
\end{align*}
and
$$(\log (g\ast h))_j =y_j+0  \text{ and } (\log (g\ast_\infty h))_j=y_j+0 \text{ for } j>2.$$
But $J_1$ is empty because $X^{\alpha}\neq 0$ only when $\alpha_2\neq 0$, for which we have assumed that $s_2=s-1$. Since $X_2\in\mathfrak{g}^{(s-1)}$ the $\vec{s}\cdot\alpha+\vec{s}\cdot\beta=s$. Consequently $(\log (g\ast h))_j=(\log (g\ast_\infty h))_j$ for every $j=1,\ldots,s$.

For every $N\in\mathbb{N}$ we set $h_1:= \exp( 2N^s X_1)$ and $h_2:=\exp (2N^{s_2}X_2)$. The discrete group $H_N:= < h_1,h_2>$ generated by those two, is isomorphic to $\mathbb{Z}^2$. We define an action on $G$: for every $(n,m)\in \mathbb{Z}^2$ and $g\in G$ we set $(n,m)\ast_N g=h_1^nh_2^m\ast g$. From the previous exposition we know that $h_1^nh_2^m\ast g =h_1^nh_2^m\ast_\infty g$, so this action we have created is independent of the choice of group structure.

We define $\text{Box}(N):=\lbrace (t_1,\ldots, t_d)\in \mathbb{R}^d : \vert t_i \vert < N^{s_i}\rbrace=\delta_N(\text{Box}(1))$ and we set $\mathcal{S}_0= \mathcal{S}^{(N)}_{0}:=\exp (\text{Box}(N))$, where $N>0$ will be fixed later. For every $(n,m)\in \mathbb{Z}^2$, we denote $\mathcal{S}_{(n,m)}:=(n,m)\ast_{N}\mathcal{S}_0$. Let $a_1$, $a_2$ be the usual generators of $\mathbb{Z}^2$, $(1,0)$ and $(0,1)$ correspondingly. We denote as $\mathcal{S}_1 = \mathcal{S}_0 \cup \mathcal{S}_{a_1}$ and $\mathcal{S}_2=\mathcal{S}_0\cup \mathcal{S}_{a_2}$. For any edge $e=\lbrace (n,m),(n,m)+a_i\rbrace$ of the Cayley graph $G(\mathbb{Z}^2,\lbrace a_1,a_2\rbrace)$, we denote $\mathcal{S}_e = (n,m)\ast_N\mathcal{S}_i$.

In the upcoming analysis of the random graph it will be important to know how many of the sets $\mathcal{S}_e$ for different edges $e$ overlap. Our aim is to bound the distance of two elements $v,u\in \mathbb{Z}^2$ for which $\mathcal{S}_v$ and $\mathcal{S}_u$ have non-trivial intersection. By obtaining such a bound, it is ensured that whenever two elements are further away, the random graphs contained in the corresponding boxes are independent. The bound that we will obtain will be independent from $N$ the size of the box.

\begin{note}
We aim to create a matching with a bond percolation on $\mathbb{Z}^2$. An edge $e$ will be open if $\mathcal{S}_e$ will contain a unique linear size connected component. This percolation will be $k$-independent, this $k$ though will depend on the number of overlaps of translations of $\mathcal{S}_e$. This is the reason why we want to bound this number independently from $N$. An additional reason we want to bound it is the following: we will prove that as $N$ grows the probability for $\mathcal{S}_e$ to contain a unique linear component also grows, we want it to be larger than the constant given by Proposition \ref{Lig}, therefore $N$ and $k$ must be independent.
\end{note}

There exists a $K>0$ independent of $N$ such that for every $v,u\in\mathbb{Z}^2$ at graph distance larger than $K$ the sets $v\ast_N \mathcal{S}_i\cap w\ast_N\mathcal{S}_j=\emptyset$, for any $i,j=0,1,2$. 

Without any loss of generality we study the case $v=(0,0)$, since every other case is a translation. Let $(n,m)\in\mathbb{Z}^2$ and $g\in (n,m)\ast_N\mathcal{S}_0$ then $g=h^n_1\ast_\infty h_2^m\ast_\infty x$ for an $x\in \mathcal{S}_0$.
We distinguish two cases.
\begin{itemize}
\item If $\vert m \vert \geq 2$, then
\begin{align*}\vert ( \log (g))_2\vert = \vert m2N^{s_2} + (\log (x))_2\vert\geq 4N^{s_2}-\vert (\log (x))_2\vert \geq 3N^{s_2},
\end{align*}
whence $g\notin \mathcal{S}_0$.
\item If $\vert m \vert <2$, then 
\begin{align*}
\vert (\log (g) )_1 \vert &= \vert n2N^s + (\log (x) )_1 + \sum_{i\in I_1} C_{\alpha,\beta}m2N^{s_2} (x)_i\vert \\
&\geq n2N^s-N^s- \sum_{i\in I_1} \vert C_{\alpha, \beta}\vert 2N^s,
\end{align*}
so that if $n>1+\sum \vert C_{\alpha,\beta}\vert$ then $\vert ( \log g)_1\vert > N^s$, hence $g\notin \mathcal{S}_0$.
\end{itemize}
Therefore, if $v$ is at graph distance larger than $K:= 2 + \sum\vert C_{\alpha,\beta}\vert $ from the $(0,0)$ then $\mathcal{S}_0\cap \mathcal{S}_v=\emptyset$.

We proceed by studying the graph inside $\mathcal{S}_i$ by defining a new finite random graph in there. 
As already mentioned, the way to study the graph locally will be through the construction of Definition \ref{2.inho}. We will define three graphs. Set $\mathcal{S}=\mathcal{S}_i$, for $i=0,1,2$, $V_\rho = \delta_{1/r}(\Gamma)\cap \mathcal{S}_i$, $\kappa^{(i)}=\lambda \mathds{1}_{B_\infty}(h^{-1}\ast_\infty g)$ and $\kappa^{(i)}_\rho = \lambda \mathds{1}_{B_\infty}(\delta_{1/r}(\delta_r(h^{-1})\ast \delta_r (g)))$. From Proposition \ref{Prop1}, we obtain inhomogeneous random graphs $G(\rho,\kappa^{(i)}_\rho)$ for $i=0,1,2$.
We connect each $g,h\in \delta_{1/r}(\Gamma)\cap \mathcal{S}_i$ with probability $\lambda \mathds{1}_{B_\infty}(\delta_{1/r}(\delta_r(h^{-1})\ast \delta_r (g)))$, giving exactly the distribution of $G'(\Gamma)$ restricted in $\mathcal{S}_i$. Thus, studying $G(p,\kappa_\rho^{(i)})$ is equivalent to studying $G'(\Gamma)$ locally.

We will bound below the norm of the integral operator $T_{\kappa^{(i)}}$ for every $i=0,1,2$ in order to apply Proposition \ref{BolProp1} and Proposition \ref{BolProp2}. For any function $f:\mathcal{S}_i\to \mathbb{R}$, the operator is
$T_{\kappa^{(i)}}f(g)=\lambda\int\limits_{\mathcal{S}_i}\mathds{1}_{B_\infty}(h^{-1}\ast g) f(h)\mathcal{H}(\mathrm{d} h)$.
Let $f\equiv 1$ then:

\begin{align*}
\Vert f\Vert_2^2 &= \int_{\mathcal{S}_i} f^2 \mathrm{d} \mathcal{H}=\mathcal{H}(\mathcal{S}_i)\\
\Vert T_{\kappa^{(i)}} f\Vert_2^2 &= \int_{\mathcal{S}_i} \left(\lambda \int_{\mathcal{S}_i} \mathds{1}_{B_\infty} (h^{-1}\ast_\infty g) \mathcal{H}(\mathrm{d} h)\right)^2 \mathcal{H}(\mathrm{d} g)=\\
&= \int_{\mathcal{S}_i} \lambda^2 \mathcal{H}(g\ast_\infty B_\infty\cap \mathcal{S}_i)^2 \mathcal{H}(\mathrm{d} g)
\end{align*}
Thus, 
\begin{align*}\Vert T_{\kappa^{(i)}}\Vert^2\geq \frac{\Vert T_{\kappa^{(i)}} f\Vert_2^2}{\Vert f\Vert_2^2}
=\lambda^2 \frac{\int_{\mathcal{S}_i} \mathcal{H}(g\ast_\infty B_\infty\cap \mathcal{S}_i)^2 \mathcal{H}(\mathrm{d} g)}{\mathcal{H}(\mathcal{S}_i)}.
\end{align*}

For $g\in \mathcal{S}_i$ with $d_{CC}(g,\partial \mathcal{S}_i)\geq 1$, it is true that $g\ast_\infty B_\infty \subset \mathcal{S}_i$, therefore $\mathcal{H}(g\ast_\infty B_\infty\cap \mathcal{S}_i)=\mathcal{H}(g\ast_\infty B_\infty)=1$. We split $\mathcal{S}_i$ into $U_i:= \lbrace g\in \mathcal{S}_i : d_{CC}(g,\partial \mathcal{S}_i)\leq 1 \rbrace$ and its complement and we obtain 
\begin{align*}
\Vert T_{\kappa^{(i)}}\Vert^2 &\geq \lambda^2 \left( \frac{\int_{U_i} \mathcal{H}(g\ast_\infty B_\infty\cap \mathcal{S}_i)^2 \mathcal{H}(\mathrm{d} g)}{\mathcal{H}(\mathcal{S}_i)}+
\frac{\int_{\mathcal{S}_i\setminus U_i}1  \mathcal{H}(\mathrm{d} g)}{\mathcal{H}(\mathcal{S}_i)}
\right)\\
&\geq \lambda^2 \frac{\mathcal{H}(\mathcal{S}_i\setminus U_i)}{\mathcal{H}(\mathcal{S}_i)}
\end{align*}

The sets $\mathcal{S}_i$ depend on the choice of $N$, we need to find an $N_0>0$ such that for any $N\geq N_0$ the operator norm $\Vert T_{\kappa^{(i)}}\Vert > 1$ for all $i=0,1,2$. Since $\lambda >1$, it is sufficient to prove that the ratio $ \mathcal{H}(\mathcal{S}_i\setminus U_i)/\mathcal{H}(\mathcal{S}_i)$ converges to $1$, as $N\to\infty$. Rather than calculating $\mathcal{H}(\mathcal{S}_i \setminus U_i)$ directly, we will find a subset $A$ of $\mathcal{S}_i \setminus U_i$ such that $\mathcal{H}(A)/\mathcal{H}(\mathcal{S}_i)\to 1$; this clearly suffices. For $\mathcal{S}_0$ this set will be $\mathcal{S}_0^{(N-L)}$ the image of $\text{Box}(N-L)$ to $G$ through the exponential. For some fixed $L$ depending on the group $G$, we will show that $\mathcal{S}_0^{(N-L)}\subset \mathcal{S}_0\setminus U$ for large $N$.

In order to prove that $\mathcal{S}_0^{(N-L)} \subseteq \mathcal{S}_0\setminus U_0$, it is sufficient to prove that if an element $g\in G$ is at distance $d_{CC}(g,\mathcal{S}_0^{(N-L)} )\leq 1$, then $g$ is contained in $\text{int}(\mathcal{S}_0)$. For any such element $g$ there exists $x\in \mathcal{S}_0^{(N-L)} $ such that $g\in B_{CC}(x,1)$, hence $g\in x\ast_\infty B_\infty$. Thus there exists a $y\in B_\infty$ for which $g=x\ast_\infty y$. Let $\log (x)=:X=(x_1,\ldots,x_d)$ and $\log (y)=:Y=(y_1,\ldots,y_d)$. Since $B_\infty$ is bounded, we know that $x_i < C^{s_i}$ for a constant $C>0$, fixed uniformly for all elements in $B_\infty$, and $y_i< (N-L)^{s_i}$. From the BCH formula (\ref{BHC3}) we obtain
\begin{align*}(X\bullet_\infty Y)_i= x_i+y_i+\sum_{I_i}C_{\alpha,\beta} X^\alpha Y^\beta\leq (N-L)^{s_i} + C^{s_i}+\sum_{I_i}C_{\alpha,\beta} (N-L)^\alpha C^\beta.
\end{align*}
We observe that $\sum_{I_i}c_{\alpha,\beta} (N-L)^\alpha C^\beta $ is a polynomial of $N$ with degree $s_i-1$ and the maximum degree is attained for $\vec{s}\beta = 1$. If $L$ is choosen such that $L>\sum_{i=1}^d\sum_{I_i}\vert C_{\alpha,\beta}\vert C$ then there is $N_0$ large enough such that for every $N>N_0$

$$
LN^{s_i-1} >\sum_{j=2}^{s_i}\left(\frac{s_i}{j}\right) L^j N^{s_i-j}(-1)^j+C^{s_i}+\sum_{I_i}c_{\alpha,\beta} (N-L)^\alpha C^\beta  
$$
we add $N^{s_i}$ and abstract $LN^{s_i-1}$ in both sides and obtain
$$
N^{s_i} > (N-L)^{s_i} + C^{s_i}+\sum_{I_i}c_{\alpha,\beta} (N-L)^\alpha C^\beta  
$$
for every $i=1,\ldots,d$. Thus $x\ast_\infty y\in \text{int}(\mathcal{S}_0)$.

Now we will show that the ratio of $ \mathcal{H}(\mathcal{S}_0^{(N-L)} )/\mathcal{H}(\mathcal{S}_0)$ tends to $1$ as $N\to\infty$. The Haar measure $\mathcal{H}$ is defined as the push-forward of the Lebesgue measure $\mu$ through $\Psi_\infty$ divided by $c_S$. It is possible to change coordinate system from canonical coordinates of second kind to exponential. Theorem \ref{LebHar} provides a relation for the corresponding measures: $\mathcal{H}(\cdot)=\vert \det \mathbf{J}_{(\Psi^{-1}_\infty \circ \exp)}(0)\vert \mu (\exp^{-1}(\cdot))$. The ratio is written as

\begin{align*}
\frac{ \mathcal{H}(\mathcal{S}_0^{(N-L)} )}{ \mathcal{H}(\mathcal{S}_0^{(N)} )}&=
\frac{\vert \det \mathbf{J}_{(\Psi^{-1}_\infty \circ \exp)}(0)\vert}{\vert \det \mathbf{J}_{(\Psi^{-1}_\infty \circ \exp)}(0)\vert}
\frac{ \mu(\exp^{-1}(\exp (\text{Box}(N-L))))}{ \mu(\exp^{-1}(\exp(\text{Box}(N))))}\\
&=\left(\frac{\mu(\text{Box}(N-L))}{\mu(\text{Box}(N))}\right)^{d_\Gamma}\\
&=\left(\frac{N-L}{N}\right)^{d_\Gamma}\to 1.
\end{align*}
This concludes the proof for $\mathcal{S}_0$.

To pass the property to $\mathcal{S}_i$ it is enough to take the translations of $\mathcal{S}_0^{(N-L)}$ through $a_i$, the generators of $\mathbb{Z}^2$. The action of $a_i$ preserves the distance $d_{CC}$, thus 
$a_i\ast_N  \mathcal{S}_0^{(N-L)} \subset \mathcal{S}_i\setminus U_i$. We proved that the measure of $\left(\mathcal{S}_0^{(N-L)}\right)^c\cap \mathcal{S}_0$ over $\mathcal{H}(\mathcal{S}_i)$ tends to $0$. Since $\mathcal{H}$ is the Haar measure the same applies for the translation $\left(a_i\ast_N\mathcal{S}_0^{(N-L)}\right)^c\cap \mathcal{S}_{a_i}$, hence the measure of the union of those two sets over $\mathcal{H}(\mathcal{S}_i)$ converges to $0$. Therefore, the measure of the complementary set in $\mathcal{S}_i$ over $\mathcal{H}(\mathcal{S}_i)$ converges to $1$ and
$$
\left( \left(\left(\mathcal{S}_0^{(N-L)}\right)^c\cap \mathcal{S}_0\right)
\cup
\left( \left(a_i\ast_N\mathcal{S}_0^{(N-L)}\right)^c\cap \mathcal{S}_{a_i} \right) \right)^c
\subset
\mathcal{S}_0^{(N-L)} \cup a_i\ast_N\mathcal{S}_0^{(N-L)}.
$$
Therefore we obtain the limit
$$\frac{\mathcal{H}\left(\mathcal{S}_0^{(N-L)} \cup a_i\ast_N \mathcal{S}_0^{(N-L)}\right)}{\mathcal{H}(\mathcal{S}_i)}\to 1$$
This conclude the proof for $\mathcal{S}_i$.

The distribution of $G(\rho,\kappa_\rho)$ is the same as $G'(\Gamma)\cap \mathcal{S}_i$. The assumptions of Theorem \ref{BolProp1} are satisfied for $(\mathcal{S}_0,\kappa^{(0)}_\rho,(V^0_\rho)_{\rho>0})$, whence there is a constant $\alpha>0$ such that $\frac{1}{\rho}C_1\left(G(\rho,\kappa^{(0)}_\rho)\right)\overset{p}\to \alpha$ as $\rho\to\infty$. Moreover the assumptions of Theorem \ref{BolProp2} are satisfied for $(\mathcal{S}_i,\kappa_\rho^{(i)},(V^i_\rho)_{\rho>0})$ for both $i=1,2$, hence $C_2\left(G(\rho,\kappa^{(i)}_\rho)\right)=o_p(\rho)$, which means that for every $\epsilon >0$ the probability $\mathbb{P}\left(C_2\left( G(\rho,\kappa^{(i)}_\rho)\right) \geq \epsilon \rho \right) \to 0$ as $\rho \to \infty$.

For every $(n,m)\in\mathbb{Z}^2$, the $G'(\Gamma)\cap \mathcal{S}_{(n,m)}$ has the same distribution as a random graph as $(n,m)\ast_N G_{\rho}^{(0)}$. Thus $\mathbb{P}(C_1
(G'(\Gamma)\cap\mathcal{S}_{(n,m)})\geq (\alpha-\epsilon)\rho)\to 1$ as $\rho\to\infty$.
Similarly for the edges, given $e=\lbrace (n,m),(n,m)+a_i\rbrace$ for $i=1,2$, the graph $G'(\Gamma)\cap \mathcal{S}_{e}$ has the same distribution as $(n,m)\ast_N (G'(\Gamma)\cap \mathcal{S}_i)$.
Therefore, for every $\epsilon>0$ we have $\mathbb{P}(C_2(G'(\Gamma)\cap \mathcal{S}_{e})\geq \epsilon\rho)\to 0$ as $\rho\to\infty$.

Now we have all the tools to establish a matching between our model and a bond percolation on $\mathbb{Z}^2$. Let $e=\lbrace (n,m),(n,m)+a_i\rbrace$ be an edge of the lattice $\mathbb{Z}^2$, and define $X(e)$ to be the event that $C_1(G'(\Gamma))\cap \mathcal{S}_{(n,m)}$ and $C_1(G'(\Gamma))\cap \mathcal{S}_{(n,m)+a_i}$ are both larger than $\alpha \rho$ while $C_2(G'(\Gamma))\cap \mathcal{S}_e < \alpha \rho$. The previous analysis proves exactly that $\mathbb{P}(X(e))\to 1$ as $r\to \infty$ and that the convergence is uniform on edges. Let $p_K$ be the constant of Proposition \ref{Lig}, and choose $r>0$ large enough such that $\mathbb{P}(X(e))\geq p_K$ for every edge $e$ of $G(\mathbb{Z}^2,\lbrace a_1,a_2\rbrace)$.

This event will provide us with the coupling between the spread-out percolation on $G$ and a $K$-independent Bernoulli percolation on $G(\mathbb{Z}^2,\lbrace a_1,a_2\rbrace)$.
We declare an edge $e$ of $G(\mathbb{Z}^2,\lbrace a_1,a_2\rbrace)$ open if $X(e)$ holds. Since for every two edges $e,f$ at distance larger than $K$ the graphs $G'(\Gamma)\cap \mathcal{S}_{e}$ and $G'(\Gamma)\cap \mathcal{S}_{f}$ have independent distributions, the events $X(e)$ and $X(f)$ are independent. As a consequence, the bond percolation defined is $K$-independent, and for any edge $e$ the probability that $e$ is open is larger than $p_K$. Since $p_K$ is chosen based on Proposition \ref{Lig}, with probability 1 there exists an infinite sequence $v_1,v_2\ldots\in \mathbb{Z}^2$ such that $v_i$ is adjacent to $v_{i+1}$ and $X(\lbrace v_i,v_{i+1}\rbrace)$ holds for every $i\in \mathbb{N}$.

For each $i\in\mathbb{N}$ the graph $G'(\Gamma) \cap \mathcal{S}_{v_i}$ contains at least one connected component of size larger than $\alpha\rho$; let $\mathcal{C}_i$ be such a component. Since $X(\lbrace v_i,v_{i+1}\rbrace)$ is true, $\mathcal{S}_{\lbrace v_i,v_{i+1}\rbrace}$ contains at most one connected component of size bigger than $\alpha\rho$, consequently $\mathcal{C}_i\cup \mathcal{C}_{i+1}$ must be a connected set.
This statement holds for every $i\in\mathbb{N}$, thus $\bigcup_{i=1}^\infty\mathcal{C}_i$ is a connected, infinite set. We conclude that, with probability 1, there is an infinite cluster in $G'(\Gamma)$. As mentioned above this graph has the same distribution as $G^{r,\lambda}(\Gamma)$ rescaled  by $\delta_{1/r}$, thus $G^{r,\lambda}(\Gamma)$ also contains an infinite connected component with probability 1.
\end{proof}

Now we provide the proofs of the lemmas.

\begin{proof}[Proof of Lemma \ref{Bounded}]
Let $x\in F_r^{-1}(\mathcal{S})$, from the previous Lemma \ref{Psilim} there is $r_1>0$ such that for $r>r_1$, $F_r(x)$ is $\eta$ $\Vert\cdot\Vert_\infty$-close to $\Psi_\infty(x)$, hence $\Psi_\infty(x)\in\mathcal{N}_\eta(\mathcal{S})$, thus $x\in \Psi^{-1}_{\infty}(\mathcal{N}_\eta(\mathcal{S}))$ which proves one part of the desired relation. 

Let $x\in\Psi^{-1}_\infty(\mathcal{S}\setminus \mathcal{N}_\eta(\partial\mathcal{S}))$, similarly from the Lemma \ref{Psilim} there is $r_2>0$ such that for every $r>r_2$ $F_r(x)$ is $\eta$ $\Vert \cdot \Vert_\infty$-close to $\Psi_\infty(x)$, hence $F_r(x)\in \mathcal{N}_\eta(\mathcal{S}\setminus\mathcal{N}_\eta(\partial\mathcal{S})) \subseteq \mathcal{S}$.

The first $r_1$ sets are bounded from its union, as a finite union of bounded sets and the rest of the sequence is contained in $\Psi^{-1}_\infty(\mathcal{N}_\eta(\mathcal{S}))$, which is a bounded set. Therefore there is a bounded set $B\subset \mathfrak{g}$ such that $F_r^{-1}(\mathcal{S})\subseteq B$ for every $r>0$.
\end{proof}

\begin{proof}[Proof of Lemma \ref{Psilim}]
First we remark that $\delta_{1/r}(\Psi_\infty (\delta_r(x)))=\Psi_\infty(x)$, combined with Lemma \ref{PsiInf} we obtain that $(\log F_r(x))_i = (\Psi_\infty(x))_i + r^{-s_i}Q_i(r^{s_{i+1}}x_{i+1},\ldots,r^{s_d}x_d)$, but \linebreak $ r^{-s_i}Q_i(r^{s_{i+1}}x_{i+1},\ldots,r^{s_d}x_d)<r^{-1}Q_i(x_{i+1},\ldots,x_d)$ since $Q_i$ has weighted degree strictly smaller than $s_i$. The number $Q_i(x_{i+1},\ldots , x_d)$ is fixed as $r$ grows to infinity, thus $r^{-s_i}Q_i(r^{s_{i+1}}x_{i+1},\ldots,r^{s_d}x_d)\to 0$ as $r\to \infty$.

If we restrict to a bounded subset $\mathcal{S}\subset G$ then there are $M_i>0$ such that
$Q_i(x_{i+1},\ldots,x_{d})<\sum_{L_i} \vert D_\alpha \vert M_{i+1}^{\alpha_{i+1}}\ldots M_d^{\alpha_d}=M$. Thus $\vert (\log F_r(x))_i - (\Psi_\infty(x))_i \vert < r^{-1}M$ for every $x\in \mathcal{S}$ and the convergence is uniform with respect to the topology induced by $\Vert \cdot \Vert_\infty$, but as already mentioned it coincides with the topology of the manifold.
\end{proof}

\begin{proof}[Proof of Lemma \ref{ExtraLe}]
There are $X_\rho, Y_\rho\in \mathfrak{g}$ such that $\exp X_\rho=g_\rho$ and $\exp Y_\rho=h_\rho$, where $X_\rho=(x^{(\rho)}_1,\ldots,x^{(\rho)}_d)$ and $Y_\rho=(y^{(\rho)}_1,\ldots,y^{(\rho)}_d)$. Since $\log$ is homeomorphism, $x^{(\rho)}_j \to x_j$ and $y^{(\rho)}_j\to y_j$ for every $j=1,\ldots,d$.

Let $Z_\rho=(z^{(\rho)}_1,z^{(\rho)}_2,\ldots,z^{(\rho)}_d)$, $W_\rho=(w^{(\rho)}_1,\ldots,w^{(\rho)}_d)$  and $W=(w_1,\ldots,w_d)$ be the elements in $\mathfrak{g}$ defined by the relations:
\begin{align*}
Z_\rho:&=\log\left(\delta_r(g_\rho)\ast \delta_r(h_\rho)\right)\\
W_\rho:&=\log\left(g_\rho\ast_\infty h_\rho\right)\\
W:&=\log(g\ast_\infty h).
\end{align*}
Since $\ast_\infty$ and $\log$ are continuous $g_\rho\ast_\infty h_\rho$ converges to $g\ast_\infty h$, hence $w^{(\rho)}_j\to w_j$ for every $j=1,\ldots, d$. Moreover
$
\delta_{1/r}(\delta_r(g_\rho)\ast\delta_r(h_\rho)) =\delta_{1/r}(\exp(Z_\rho))=\exp (\delta_{1/r} (Z_\rho))
$ and $\delta_{1/r} (Z_\rho)=(r^{-s}z_1,\ldots,r^{-1}z_d)$. From BCH formula (\ref{BHC2}) on $\delta_r(X_\rho)$ and $\delta_r(Y_\rho)$ and the BCH formula (\ref{BHC3}) on $X_\rho$ and $Y_\rho$ we obtain that
\begin{align*}
z^{(\rho)}_j &=r^{s_j}x_j^{(\rho)} + r^{s_j}y_j^{(\rho)}+ \sum_{I_j}C_{\alpha,\beta} r^{\alpha} X_\rho^{\alpha}r^{\beta}Y_\rho^{\beta}+\sum_{J_j}D_{\alpha,\beta} r^{\alpha} X_\rho^{\alpha}r^{\beta}Y_\rho^{\beta}\\
&=r^{s_j }w^{(\rho)}_j+\sum_{J_j}D_{\alpha,\beta} r^{\alpha+\beta} X_\rho^{\alpha} Y_\rho^{\beta}.
\end{align*}
Therefore $r^{-s_j}z_j^{(\rho)}-w_j^{(\rho)}\to 0$, because $\sum_{J_j}D_{\alpha,\beta} r^{\alpha+\beta-s_j} X_\rho^{\alpha} Y_\rho^{\beta} \to 0$. We conclude that $\delta_{1/r} (Z_\rho)\to W$, as $\rho \to \infty$, recall that $\rho=c_S r^{d_\Gamma}$. Since $\exp$ is a homeomorphism $\exp (\delta_{1/r} Z_\rho)\to \exp W$, but $\exp (\delta_{1/r} (Z_\rho))$ is exactly $\delta_{1/r}(\delta(g_\rho)\ast\delta_r(h_\rho))$. 
\end{proof}

\section{Transitive graphs}\label{6}

As we mentioned in the introduction, in proving \cref{transitive}, we make use of a theorem of Trofimov \cite{trofimov} stating that an arbitrary transitive graph of polynomial growth can be approximated by a Cayley graph in a certain precise sense. There are various ways of formulating this result, but for us the most convenient will be the version appearing in work of Tessera and the second author \cite[Theorem 2.1]{tt.trof}, which essentially follows from an argument of Woess \cite{woess}. Before we state this theorem, we present some background. This is mostly reproduced from \cite{tt.trof}, where the reader can find full proofs; much of that is in turn based on \cite{woess}.

If $G$ is a transitive graph and $H<\Aut(G)$ is a subgroup then we define $G/H$ to be the quotient graph with vertices $\{H(x):x\in G\}$, and $H(x)\sim H(y)$ in $G/H$ if and only if there exist $x_0\in H(x)$ and $y_0\in H(y)$ such that $x_0\sim y_0$ in $G$. If $\Gamma<\Aut(G)$ is a transitive subgroup that normalizes $H$ then the quotient graph $G/H$ is invariant under the action of $\Gamma$ on $G$, and the action of $\Gamma$ on $G$ descends to an action of $\Gamma$ on $G/H$. We write $\Gamma_{G/H}$ for the image of $\Gamma$ in $\Aut(G/H)$ induced by this action; thus $\Gamma_{G/H}$ is the quotient of $\Gamma$ by the normal subgroup $\{g\in\Gamma:gH(x)=H(x)\text{ for every }x\in G\}$.

The automorphism group $\Aut(G)$ is a metrizable topological group under the topology of pointwise convergence, in which vertex stabilizers are compact. Given a subgroup $\Gamma<\Aut(G)$ and a vertex $x\in G$, we write $\Gamma_x$ for the stabilizer of $x$. The vertex stabilizers $\Gamma_x$ are all isometric to one another, and in particularly have the same diameter and cardinality; $\Gamma$ is discrete if and only if this cardinality is finite.

A subset $U\subset\Aut(G)$ is precompact if and only if all its orbits are finite, if and only if at least one orbit is finite. In particular, if $H<\Aut(G)$ is a compact subgroup normalized by a transitive subgroup $\Gamma<\Aut(G)$ then the fibres of the projection $G\to G/H$ are finite.

If $\Gamma$ acts transitively on $G$ then, given a fixed vertex $o\in G$, the set $S=\{g\in\Gamma:d(g(o),o)\le1\}$ is a compact symmetric generating set for $\Gamma$. The fact that $S$ contains the identity has the notationally convenient effect of ensuring that for each $n\in\N$ the ball of radius $n$ in $G(\Gamma,S)$ is exactly the $n$-fold product set $S^n=\{s_1\cdots s_n:s_i\in S\}$. Moreover, given $g,h\in\Gamma$ and a non-negative integer $n$, we have $d(g(o),h(o))\le n$ if and only if there exists $s\in S^n$ such that $g=hs$, and given a (left-invariant) Haar measure $\mu$ on $\Gamma$ and a vertex $x\in G$, we have
\begin{equation}\label{eq:VTgrowth}
\beta_G(n)=\frac{\mu(S^n)}{\mu(\Gamma_x)}.
\end{equation}

\begin{lemma}\label{image.of.stab}
Let $G$ be a transitive graph, let $H\lhd\Gamma<\Aut(G)$, and write $\pi:\Gamma\to\Gamma_{G/H}$ for the quotient homomorphism. Let $o\in G$. Then $(\Gamma_{G/H})_{H(o)}=\pi(\Gamma_o)$.
\end{lemma}
\begin{proof}
This is claimed in \cite[Lemma 3.5]{tt.trof}, but it does not quite appear to follow immediately from the preceding conclusion, which immediately implies only that $\pi(g)\in(\Gamma_{G/H})_{H(o)}$ if and only if there exists $h\in H$ such that $hg\in\Gamma_o$. However, \cite[Lemma 3.6]{tt.trof} then implies that this occurs if and only if there exists $h\in\ker\pi$ such that $hg\in\Gamma_o$, as required.
\end{proof}

We will actually need the following slight refinement of \cite[Theorem 2.1]{tt.trof}, in which the conclusion that the action of the finite-index nilpotent subgroup on the quotient graph is free replaces the weaker condition that $\Gamma_{G/H}$ acts with finite vertex stabilizers.
\begin{theorem}\label{thm:trof}
Let $G$ be a connected, locally finite transitive graph of polynomial growth, and let $o\in G$. Let $\Gamma<\Aut(G)$ be a transitive subgroup. Then there is a compact normal subgroup $H\lhd\Gamma$ such that $\Gamma_{G/H}$ contains a normal torsion-free nilpotent subgroup of finite index that acts freely on $G/H$.
\end{theorem}
\begin{Remark}
By following the same argument and keeping track of all the quantitative aspects, one can also deduce that the finite-index nilpotent subgroup acts freely on the quotient graph in Tessera and the second author's finitary version of Trofimov's theorem \cite[Theorem 2.3]{tt.trof}.
\end{Remark}

To obtain the additional conclusion that the finite-index normal nilpotent subgroup acts freely on $G/H$, we use the following lemma.
\begin{proposition}\label{nilp.free}
Let $G$ be a connected, locally finite transitive graph of polynomial growth, and let $o\in G$. Let $\Gamma<\Aut(G)$ be a transitive subgroup acting on $G$ with finite vertex stabilizers, and containing a normal nilpotent subgroup $N$ with finite index. Then there exists a finite normal subgroup $H\lhd\Gamma$ such that the image of $N$ in $\Gamma_{G/H}$ acts freely on $G/H$.
\end{proposition}
\begin{proof}
If $N_o=\{\id\}$ then $N$ acts freely on $G$, since if $n(v)=v$ for some $n\in N$ and $v\in G$ then by transitivity there exists $g\in \Gamma$ such that $g(v)=o$, so that $gng^{-1}(o)=o$ and hence $n=\id$ by normality of $N$. We may therefore assume that $N_o$ contains at least one non-identity element. In particular, we may also assume that $|\Gamma_o|>1$ and proceed by induction on $|\Gamma_o|$.

Consider $N_o^\Gamma$, the normal closure in $\Gamma$ of the stabilizer $N_o$. Let $S$ be the symmetric generating set $\{g\in\Gamma:d_G(g(o),o)\le1\}$, which is a finite symmetric generating set for $\Gamma$. Since $S$ is finite, $S^2$ is an \emph{approximate group}, which simply means that $S^2$ is symmetric and contains the identity and that there is a finite set $X$ such that $S^2\subset XS$. We may therefore apply \cite[Proposition 6.5]{tt.trof} to conclude that $N_o^\Gamma$ is finite.

Write $\pi:\Gamma\to\Gamma_{G/N_o^\Gamma}$ for the quotient homomorphism, noting that $(\Gamma_{G/N_o^\Gamma})_{N_o^\Gamma(o)}=\pi(\Gamma_o)$ by \cref{image.of.stab}, and hence $|(\Gamma_{G/N_o^\Gamma})_{N_o^\Gamma(o)}|<|G_o|$ because $N_o$ contains a non-trivial element. The proposition therefore follows from applying the induction hypothesis to the group $\Gamma_{G/N_o^\Gamma}$ acting on the graph $G/N_o^\Gamma$.
\end{proof}
The conclusion that the finite-index normal nilpotent subgroup is torsion-free is automatic by the following standard lemma.
\begin{lemma}\label{virt.nilp->normal}
Suppose $\Gamma$ is a finitely generated group and $N$ is a nilpotent subgroup of finite index. Then $N$ contains a torsion-free subgroup of finite index that is normal in $\Gamma$.
\end{lemma}
\begin{proof}
A subgroup of finite index in a finitely generated group is finitely generated, and every finitely generated nilpotent group embeds as a finite index subgroup in $A\times B$ for some finite group $A$ and some torsion-free group $B$ (see \cite[Theorem~2.1]{Bau}), so $N$ contains a further finite-index subgroup $H$ that is torsion-free. It is easy to check that the subgroup $\bigcap_{gH\in\Gamma/H}gHg^{-1}$ is a well-defined finite-index normal subgroup of~$\Gamma$.
\end{proof}

\begin{proof}[Proof of \cref{thm:trof}]
It is stated in \cite[Theorem 2.1]{tt.trof} that there is a compact normal subgroup $H_0\lhd\Gamma$ such that $\Gamma_{G/H_0}$ is virtually nilpotent and acts on $G/H_0$ with finite vertex stabilizers. \cref{virt.nilp->normal} implies that $\Gamma$ contains a normal nilpotent subgroup $N$ of finite index. \cref{nilp.free} therefore implies that there exists a normal subgroup $H\lhd\Gamma$ containing $H_0$ such that $[H:H_0]<\infty$, so that $H$ is compact, and such that the image of $N$ in $\Gamma_{G/H}$ acts freely on $G/H$. The image of $N$ is a nilpotent subgroup with finite index in $\Gamma_{G/H}$, so by \cref{virt.nilp->normal} it contains a torsion-free nilpotent subgroup with finite index that is normal in $\Gamma_{G/H}$.
\end{proof}

Our aim in applying \cref{thm:trof} is to compare spread-out percolation on $G$ to spread-out percolation on $G/H$, and spread-out percolation on $G/H$ to spread-out percolation on $\Gamma_{G/H}$. We start by comparing spread-out percolation on $G$ to spread-out percolation on $G/H$.

Benjamini and Schramm \cite[Theorem 1]{bperc96} famously showed that for each $v\in G$ and $p\in[0,1]$ the image under $\pi$ of the cluster of $v$ in Bernoulli-$p$ percolation on $G$ stocastically dominates the cluster of $\pi(v)$ in Bernoulli-$p$ percolation on $G / H$. However, this is of limited use in our context, since spread-out percolation with parameter $\lambda$ at scale $r$ corresponds to very different values of $p$ in Bernoulli percolation on $G/H$ or on $G$. Fortunately, however, a similar argument to Benjamini and Schramm's yields the following result. Here and from now on, we write $\Prob_p^G$ for the law of Bernoulli-$p$ bond percolation on $G$. Given a vertex $v$ of $G$ we write $K_v$ for the cluster containing $v$, and given another vertex $u$ we write $u\leftrightarrow v$ for the event that $u$ and $v$ belong to the same open cluster.

\begin{proposition}\label{Benj-Schr}
Let $G$ be a locally finite graph, let $H<\Aut(G)$, and let $\pi:G\to G/H$ be the quotient map. Let $r\in\N$, and suppose that each orbit of $H$ has size $k\in\N$ and diameter at most $\ell\in\N$. Then for each $v\in G$ and $p\in [0,1]$, the image under $\pi$ of the cluster of $v$ in Bernoulli-$p$ percolation on $G_{r+\ell}$ stocastically dominates the cluster of $\pi(v)$ in Bernoulli-$(1-(1-p)^k)$ percolation on $(G / H)_r$. That is,
\[
\Prob_p^{G_{r+\ell}}\Bigl(\pi(K_v) \in A\Bigr)\ge\Prob_{1-(1-p)^k}^{(G/H)_r}\Bigl( K_{\pi(v)} \in A\Bigr)
\]
for every increasing measurable set $A \subseteq \{0,1\}^{G/H}$.
\end{proposition}
\begin{proof}
We describe a procedure for coupling Bernoulli-$(1-(1-p)^k)$ bond percolation on $(G/H)_r$ to Bernoulli-$p$ bond percolation on $G_{r+\ell}$, and exploring the cluster of $\pi(v)$ in $(G/H)_r^{1-(1-p)^k}$ and the cluster of $v$ in $G^{(p)}_{r+\ell}$. We define random sequences $C_0,C_1,C_2,\ldots$ of subsets of $G/H$ and $B_0,B_1,B_2,\ldots$ of sets of edges in $(G/H)_r$, and, for each $x\in\bigcup_iC_i$, a vertex $u_{x}$ of $G$ belonging to the orbit $x$. These sets and elements will depend on the random graph $G^{(p)}_{r+\ell}$, and will be defined such that $\bigcup_iC_i$ is the cluster of $\pi(v)$ in the coupled $(G/H)_r^{1-(1-p)^k}$, and such that for each $x\in\bigcup_iC_i$ we have $u_x\leftrightarrow v$ in $G^{(p)}_{r+\ell}$, proving the proposition.

Start with $C_0=\{\pi(v)\}$ and $B_0=\varnothing$, and set $u_{\pi(v)}=v$. Then, once $C_n$ and $B_n$ are defined, if $B_n$ contains the whole edge boundary $\partial C_n$ of $C_n$ in $(G/H)_r$ then we stop. If $\partial C_n\not\subset B_n$, let $x\in C_n$ and $y\in(G/H)\setminus C_n$ be such that $(x,y)\in\partial C_n\setminus B_n$, with $x$ chosen to be at minimum possible distance from $\pi(v)$. Now note that for each vertex $z$ of $G$ belonging to the $H$-orbit $y$ we have $d(u_x,z)\le r+\ell$, so that the graph $G_{r+\ell}$ contains the edge $(u_x,z)$. If all $k$ of these edges are closed in $G^{(p)}_{r+\ell}$, we declare $(x,y)$ to be closed in $(G/H)_r$, and set $C_{n+1}=C_n$ and $B_{n+1}=B_n\cup\{(x,y)\}$. If at least one of these edges is open, however, then we declare $(x,y)$ to be open in $(G/H)_r$, set $B_{n+1}=B_n$ and $C_{n+1}'=C_n'\cup\{y\}$, and choose $u_y\in y$ arbitrarily so that $(u_x,u_y)$ is open.

Note that each edge of $(G/H)_r$ that we explore in this way has probability $1-(1-p)^k$ of being open, and that decisions about distinct edges are independent. If we now declare edges of $(G/H)_r$ we did not explore to be open independently with probability $1-(1-p)^k$, we recover Bernoulli-$(1-(1-p)^k)$ bond percolation on $(G/H)_r$, with $K_{\pi(v)}=\bigcup_iC_i$ as required (note that by choosing $x$ to have the smallest possible distance from $\pi(v)$ we ensured that $\bigcup_iC_i$ contained every vertex of $K_{\pi(v)}$).
\end{proof}

Next, we want to compare spread-out percolation on $G/H$ to spread-out percolation on $\Gamma_{G/H}$. For this, we use the following proposition, the proof of which is somewhat similar to that of \cref{Benj-Schr}.

\begin{proposition}\label{fin.ind.action}
Let $G$ be a connected transitive graph with distinguished vertex $o$, and suppose $\Gamma<\Aut(G)$ acts transitively on $G$ and contains a 
subgroup $H<\Gamma$ that acts freely on $G$ with exactly $k\in\N$ distinct orbits. Let $S=\{g\in\Gamma:d(o,g(o))\le1\}$, let $X\subset S^k$ be a minimal set containing the identity and satisfying $G=HX(o)$, and let $\eta:G\to H$ and $\xi:G\to X$ be the unique maps such that $v=\eta(v)\xi(v)(o)$ for all $v\in G$. Then 
for each $p\in[0,1]$, $m,r\in\N$ and $h\in H$, the image under $\eta$ of the cluster of $h(o)$ in Bernoulli-$p$ bond percolation on $G_{rm+2k}$ stochastically dominates the cluster of $h$ in Bernoulli-$(1-(1-p)^k)$ bond percolation on $G(H,(H\cap S^m)^r)$. That is,
\[
\Prob_p^{G_{rm+2k}}\Bigl(\eta(K_{h(o)})\in A\Bigr)\ge\Prob_{1-(1-p)^k}^{G(H,(H\cap S^m)^r)}\Bigl(K_h\in A\Bigr)
\]
for every increasing measurable set $A\subset \{0,1\}^H$.
\end{proposition}

\begin{proof}
We describe a procedure for coupling Bernoulli-$(1-(1-p)^k)$ bond percolation on $G(H,(H\cap S^m)^r)$ to Bernoulli-$p$ bond percolation on $G_{rm+2k}$, and exploring the cluster of $h$ in $G^{(1-(1-p)^k)}(H,(H\cap S^m)^r)$ and the cluster of $h(o)$ in $G^{(p)}_{rm+2k}$. We define random sequences $C_0,C_1,C_2,\ldots$ of subsets of $H$ and $B_0,B_1,B_2,\ldots$ of sets of edges in $G(H,(H\cap S^m)^r)$, and, for each $u\in\bigcup_iC_i$, an element $x_u\in X$. These sets and elements will depend on the random graph $G^{(p)}_{rm+2k}$, and will be defined such that $\bigcup_iC_i$ is the cluster of $h$ in the coupled $G^{(1-(1-p)^k)}(H,(H\cap S^m)^r)$, and such that for each $u\in\bigcup_iC_i$ we have $ux_u(o)\leftrightarrow h(o)$ in $G^{(p)}_{rm+2k}$, proving the proposition.

Start with $C_0=\{h\}$ and $B_0=\varnothing$, and set $x_h=\id$. Then, once $C_n$ and $B_n$ are defined, if $B_n$ contains the whole edge boundary $\partial C_n$ of $C_n$ in $G(H,(H\cap S^m)^r)$ then we stop. If $\partial C_n\not\subset B_n$, let $u\in C_n$ and $s\in(H\cap S^m)^r$ be such that $(u,us)\in\partial C_n\setminus B_n$, with $u$ chosen to be at minimum possible distance from $h$. Now note that for each $y\in X$ we have $d_G(ux_u(o),usy(o))\le rm+2k$, so that the graph $G_{rm+2k}$ contains the edge $(ux_u(o),usy(o))$. If all $k$ of these edges are closed in $G^{(p)}_{rm+2k}$, we declare $(u,us)$ to be closed in $G(H,(H\cap S^m)^r)$, and set $C_{n+1}=C_n$ and $B_{n+1}=B_n\cup\{(u,us)\}$. If at least one of these edges is open, however, then we declare $(u,us)$ to be open in $G(H,(H\cap S^m)^r)$, set $B_{n+1}=B_n$ and $C_{n+1}=C_n\cup\{us\}$, and choose $x_{us}\in X$ arbitrarily so that $(ux_u(o),usx_{us}(o))$ is open.

Note that each edge of $G(H,(H\cap S^m)^r)$ that we explore in this way has probability $1-(1-p)^k$ of being open, and that decisions about distinct edges are independent. If we now declare edges of $G(H,(H\cap S^m)^r)$ we did not explore to be open independently with probability $1-(1-p)^k$, we recover Bernoulli-$(1-(1-p)^k)$ bond percolation on $G(H,(H\cap S^m)^r)$, with $K_h=\bigcup_iC_i$ as required.
\end{proof}

The assumption in \cref{fin.ind.action} that $X\subset S^k$ does not render the proposition vacuous, thanks to the following lemma.
\begin{lemma}\label{balls-orbits}
Let $G$ be a connected transitive graph with distinguished vertex $o$, and suppose $\Gamma<\Aut(G)$ acts transitively on $G$ and contains a 
subgroup $H<\Gamma$ that acts on $G$ with at least $k\in\N$ distinct orbits. Let $S=\{g\in\Gamma:d(o,g(o))\le1\}$. Then $S^{k-1}(o)$ has non-empty intersection with at least $k$ distinct $H$-orbits.
\end{lemma}
\begin{proof}
The case $k=1$ is trivial, so we may assume that $k\ge2$ and, by induction, that $S^{k-2}(o)$ has non-empty intersection with at least $k-1$ distinct $H$-orbits. Writing $U$ for the union of those orbits of $H$ having non-empty intersection with $S^{k-2}(o)$, we may also assume that $U\ne G$, so that there exists a vertex $u\in U$ with a neighbour $v\in G\setminus U$. By definition, there exists $h\in H$ such that $h(u)\in S^{k-2}(o)$, and hence $h(v)\in S^{k-1}(o)$, so that $S^{k-1}(o)$ meets at least one more $H$-orbit than $S^{k-2}(o)$ does.
\end{proof}

In order to use \cref{fin.ind.action} to compare spread-out percolation on $G/H$ to spread-out percolation on the torsion-free nilpotent subgroup of $\Gamma_{G/H}$, we need to be able to compare the sizes of balls in the two spaces. For this, we use the following lemma.

\begin{lemma}\label{HcapSlarge}
Let $\Gamma$ be a group with finite symmetric generator set $S$ containing the identity, let $H\le\Gamma$ be a finite-index subgroup, and let $X\subset S^k$ be a set of right-coset representatives for $H$ in $\Gamma$. Let $m>2k$. Then for every $r\in\N$ we have $S^{r(m-2k)}\subset (H\cap S^m)^rX$. In particular, $H\cap S^m$ generates $H$ and $|S^{r(m-2k)}|\le[G:H]\,|(H\cap S^m)^r|$.
\end{lemma}
\begin{proof}
The fact that $X$ is a set of coset representatives implies that $S^{m-k}\subset HX$. Since $X\subset S^k$ this in fact implies that
\begin{equation}\label{eq:HcapS}
S^{m-k}\subset(H\cap S^m)X,
\end{equation}
and proves the case $r=1$ a fortiori. For $r>1$, we have
\begin{align*}
S^{r(m-2k)}&=S^{(r-1)(m-2k)}S^{m-2k}\\
     &\subset (H\cap S^m)^{r-1}XS^{m-2k}&&\text{(by induction)}\\
     &\subset (H\cap S^m)^{r-1}S^{m-k}\\
     &\subset (H\cap S^m)^rX&&\text{(by \eqref{eq:HcapS})},
\end{align*}
as required.
\end{proof}

We also need to compare balls of slightly different radii in $G$ and in $G/H$, for which we use the following result of Tessera, which actually holds in the more general setting of measured metric spaces with \emph{property (M)} (see \cite{tessera.folner} for definitions and details).
\begin{proposition}[Tessera {\cite[Theorem 4]{tessera.folner}}\footnote{The statement of \cite[Theorem 4]{tessera.folner} does not include the fact that $C$ and $\delta$ depend only on $K$, but this follows from the proof.}]\label{Folner}
Given $K>1$ there exist $C=C(K)\ge1$ and $\delta=\delta(K)>0$ such that if $G$ is a locally finite transitive graph and $|\beta_G(2n)|\le K|\beta_G(n)|$ for all $n\in\N$ then
\[
|B(x,n+1)\setminus B(x,n)|\le Cn^{-\delta}\beta_G(n)
\]
for all $x\in G$ and $n\in\N$.
\end{proposition}
It will be convenient to record the following particular consequence of \cref{Folner}.
\begin{lemma}\label{Folner.cor}
Let $k,t,K\in\N$ and $\eps>0$. Then there exists $M=M(k,t,K,\eps)\in\N$ such that if $G$ is a locally finite transitive graph satisfying $|\beta_G(2n)|\le K|\beta_G(n)|$ for all $n\in\N$ then $|\beta_G(r(m-t))|\ge(1-\eps)|\beta_G((r+1)m+k)|$ for all integers $m,r\ge M$.
\end{lemma}
\begin{proof}
Applying \cref{Folner} with $G_r$ in place of $G$, provided $m$ is large enough in terms of $t$, $K$ and $\eps$ we have $\beta_G(r(m-t))\ge(1-\eps)^{1/3}\beta_G(rm)$. Applying \cref{Folner} with $G_m$ now in place of $G$, provided $r$ is large enough in terms of $K$ and $\eps$ we have $\beta_G(rm)\ge(1-\eps)^{1/3}\beta_G((r+1)m)$. Finally, provided $r$ is large enough in terms of $t$, $K$ and $\eps$, \cref{Folner} implies that $\beta_G((r+1)m)\ge(1-\eps)^{1/3}\beta_G((r+1)m+k)$.
\end{proof}

\begin{proof}[Proof of \cref{transitive}]
Let $H\lhd\Gamma$ be the compact subgroup given by \cref{thm:trof}, let $N\lhd\Gamma_{G/H}$ be the normal torsion-free nilpotent subgroup of finite index in $\Gamma_{G/H}$ acting freely on $G/H$ given by the same theorem. Since $H$ is compact, its orbits are finite, so that $G/H$ also has superlinear growth. Since $\Gamma_{G/H}$ acts with finite vertex stabilizers, the set $S=\{g\in\Gamma/H:d(H(o),g(H(o)))\le1\}$ is a finite symmetric generating set, and $\Gamma_{G/H}$ is a finitely generated group of superlinear polynomial growth by \eqref{eq:VTgrowth}. Note therefore that $N$ is not cyclic, since otherwise $N$ and hence $\Gamma_{G/H}$ would have linear growth.

Write $k\in\N$ for the number of distinct orbits of the action of $N$ on $G/H$, and let $X\subset S^k$ be a minimal set containing the identity and satisfying $G/H=NX(H(o))$, noting that such a set exists by \cref{balls-orbits}. Let $Y$ be a set of right-coset representatives for $N$ in $\Gamma/H$ such that $X\subset Y$ and such that for each $Y(H(o))=X(H(o))$ (such a set exists because, given a right coset $Ng$ of $N$ in $\Gamma$, we have $g(H(o))=nx(H(o))$ for some $n\in N$ and $x\in X$, so that $Ng=N(n^{-1}g)$ and $n^{-1}g(H(o))=x(H(o))$). Fix $t\in\N$ be such that $Y\subset S^t$.

Fix $\lambda>1$ and $m,n\in\N$ with $m>2t$, and let $q=\lambda^{3/4}/|B(H(o),(n+1)m+2k)|$. Provided $m$ and $n$ are large enough, \cref{Folner.cor} implies that $\beta_{G/H}(n(m-2t))\ge\lambda^{-1/4}\beta_{G/H}((n+1)m+2k)$, and hence that $q\ge\lambda^{1/2}/\beta_{G/H}(n(m-2t))$. It follows from \cref{HcapSlarge} that $B(H(o),n(m-2t))=S^{n(m-2t)}(H(o))\subset(N\cap S^m)^nY(H(o))=(N\cap S^m)^nX(H(o))$, hence $\beta_{G/H}(n(m-2t))\le k|(N\cap S^m)^n|$, and hence
\[
q\ge\frac{\lambda^{1/2}}{k|(N\cap S^m)^n|}.
\]
Provided $n$ is large enough, we also have $1-(1-q)^k\ge\lambda^{-1/4}kq$, and hence therefore
\[
1-(1-q)^k\ge\frac{\lambda^{1/4}}{|(N\cap S^m)^n|}.
\]
Provided $n$ is large enough, \cref{HcapSlarge,4.1} then imply that
\[
\Prob_{1-(1-q)^k}(\text{$G(N,(N\cap S^m)^n)$ contains an infinite cluster})=1,
\]
and then \cref{fin.ind.action} implies that
\[
\Prob_q(\text{$(G/H)_{nm+2k}$ contains an infinite cluster})=1.
\]
Since $(G/H)_{nm+2k}$ is a subgraph of $(G/H)_r$ for all $r\ge nm+2k$, this then implies that $(G/H)_{r,\lambda^{3/4}}$ has an infinite component with probability 1 for every integer $r\in[nm+2k,(n+1)m+2k]$. Leaving $m$ fixed but large enough, and allowing $n$ to vary over all large enough integers, this proves that $(G/H)_{r,\lambda^{3/4}}$ has an infinite component with probability 1 for every integer $r\ge m+2k$.

Now write $\ell$ for the diameter of the orbits of $H$, fix an integer $r>\ell$, and set $p=\lambda/\beta_G(r)$. For large enough $r$ we have $1-(1-p)^{|H(o)|}\ge\lambda^{-1/8}|H(o)|p$, hence
\[
1-(1-p)^{|H(o)|}\ge\frac{\lambda^{7/8}|H(o)|}{\beta_G(r)}\ge\frac{\lambda^{7/8}}{\beta_{G/H}(r)}.
\]
Moreover, \cref{Folner} implies that for large enough $r$ we have $\beta_{G/H}(r)\le\lambda^{1/8}\beta_{G/H}(r-\ell)$, and hence $1-(1-p)^{|H(o)|}\ge\lambda^{3/4}/\beta_{G/H}(r-\ell)$. Provided $r$ is large enough, the previous paragraph then imples that $\Prob_{1-(1-p)^{|H(o)|}}(\text{$(G/H)_{r-\ell}$ contains an infinite cluster})=1$,
and then \cref{Benj-Schr} implies that $\Prob_p(\text{$G_r$ contains an infinite cluster})=1$, as required.
\end{proof}

We close this section and the paper by returning to one of the first assertions we made in the introduction: that equality holds in \eqref{eq:pc>} for $G$ transitive if and only if $G$ is a tree. This is well known in folklore, but we have been unable to find an explicit reference, so we record a proof here for completeness. One can also prove this more elementarily, but we opt for brevity.
\begin{proposition}\label{prop:tree}
Suppose $G=(V,E)$ is a connected transitive graph of degree $d\ge2$. Then $p_c(G)\ge1/(d-1)$, with equality if and only if $G$ is a tree.
\end{proposition}
\begin{proof}
Equality in the case of a tree is well known (see e.g. \cite[Theorem 1.8]{lyons-peres}). To see that strict inequality holds for transitive graphs that are not trees, we use a result of Martineau and Severo \cite[Theorem 0.1]{mart-sev.strict.monotonicity}, a special case of which states that if $G_0=(V_0,E_0)$ is a transitive graph with $p_c(G)<1$ and $H$ is a non-trivial group acting freely on $V_0$ by automorphisms then $p_c(G_0/H)>p_c(G_0)$. Now simply note that the universal cover of $G$ (indeed any $d$-regular graph) is the $d$-regular tree $T_d$, and that $G$ is then the quotient of $T_d$ by the group of deck transformations, which acts freely.
\end{proof}

\footnotesize{
\bibliographystyle{abbrv}
\bibliography{vivlio3}
}

\end{document}